\documentclass[11pt]{amsart}
\usepackage{amscd,amssymb}
\usepackage[mathscr]{eucal}
\usepackage[all]{xy}

\title{twisted jets, motivic measure and \linebreak orbifold cohomology}
\author{Takehiko Yasuda }
\address{Department of Mathematical Sciences, University of Tokyo,
Komaba,  Meguro, Tokyo, 153-8914, Japan }
\email{t-yasuda@ms.u-tokyo.ac.jp}


\theoremstyle{plain}
\newtheorem{thm}{Theorem}[section]

\newtheorem{prop}[thm]{Proposition}
\newtheorem{cor}[thm]{Corollary}
\newtheorem{lem}[thm]{Lemma}

\theoremstyle{definition}
\newtheorem{defn}[thm]{Definition}
\newtheorem{expl}[thm]{Example}

\theoremstyle{remark}
\newtheorem{rem}[thm]{Remark}

\newcommand{\A}{\mathbb A}
\newcommand{\CC}{\mathbb C}
\newcommand{\QQ}{\mathbb Q}
\newcommand{\ZZ}{\mathbb Z}
\newcommand{\LL}{\mathbb L}
\newcommand{\Zz}{\mathbb Z _{\geq 0}}

\newcommand{\I}{\mathscr I}
\newcommand{\J}{\mathscr J}
\newcommand{\X}{\mathscr X}
\newcommand{\Y}{\mathscr Y}

\newcommand{\D}{\mathscr D}
\newcommand{\U}{\mathscr U}
\newcommand{\cL}{\mathscr L}
\newcommand{\cO}{\mathscr O}
\newcommand{\cK}{\mathscr K}
\newcommand{\cA}{\mathscr A}
\newcommand{\cN}{\mathscr N}
\newcommand{\cW}{\mathscr W}
\newcommand{\cF}{\mathscr F}
\newcommand{\cV}{\mathscr V}

\newcommand{\fm}{\mathfrak m}
\newcommand{\fp}{\mathfrak p}

\newcommand{\MM}{\mathbb M}

\newcommand{\Mhat}{\hat{\mathbb M}}
\newcommand{\GHS}{K_0(\mathrm{HS})}

\newcommand{\HS}{\mathrm{HS}}

\newcommand{\muGor}{\mu^{\mathrm{Gor}}}
\newcommand{\sch}{(\mathrm{Sch}/\CC)}
\newcommand{\schS}{(\mathrm{Sch}/S)}
\newcommand{\ob}{\mathrm{ob}\,}
\newcommand{\Ker}{\mathrm{Ker}\,}
\newcommand{\ord}{\mathrm{ord}\,}
\newcommand{\Hom}{\mathrm{Hom}}
\newcommand{\Spec}{\mathrm{Spec}\,}
\newcommand{\cSpec}{\mathcal{S} pec\,}

\newcommand{\GL}{\mathrm{GL}}
\newcommand{\SL}{\mathrm{SL}}
\newcommand{\Img}{\mathrm{Im}\,}

\newcommand{\Conj}{\mathrm{Conj}}

\newcommand{\Gr}{\mathrm{Gr}}
\newcommand{\Aut}{\mathrm{Aut}}
\newcommand{\pr}{\mathrm{pr}}

\newcommand{\Horb}{H_{\mathrm{orb}}}
\newcommand{\id}{\mathrm{id}}
\newcommand{\codim}{\mathrm{codim}\,}
\newcommand{\diag}{\mathrm{diag}}
\newcommand{\nor}{\mathrm{nor}}
\newcommand{\sing}{\mathrm{sing}}

\numberwithin{equation}{section}

\begin{document}

\begin{abstract}
We introduce the notion of twisted jets.
For a Deligne-Mumford stack $\X$ of finite type over $\CC$,
a twisted $\infty$-jet on $\X$ is a representable morphism 
$\D \to \X$ such that $\D$ is a smooth Deligne-Mumford stack with
 the coarse moduli space $\Spec \CC[[t]]$.
We study the motivic measure on the space of 
the twisted $\infty$-jets on a smooth Deligne-Mumford stack.

As an application, we prove that two birational minimal models
with Gorenstein quotient singularities have the same orbifold cohomology
with Hodge structure.
\end{abstract}

\maketitle

\section{Introduction}

In 1995, Kontsevich produced the theory named \textit{motivic integration} 
\cite{Orsay}. 
Since then, this remarkable idea has become a powerful 
method for examining both the local and global structures of varieties.

Let $X$ be a variety over $\CC$. 
For $ n \in \Zz \cup \{ \infty \} $, 
an $n$-\textit{jet} on $X$ is
a $\CC[[t]]/(t^{n+1})$-point of $X$, where we have followed the 
convention $(t ^\infty ) = (0) $.
The  $n$-jets of $X$ naturally constitute a variety 
(or pro-variety if $n = \infty $), denoted $L_n X$.
For $ m \ge n$,  
the natural surjection $\CC[[t]]/(t^{m+1}) \to \CC[[t]]/(t^{n+1})$
induces the truncation morphism $ L_m X \to L_n X $.

Consider the case where $X$ is smooth and of dimension $d$. 
Then $L_n X$ is a locally trivial affine space bundle 
over $X$. (Whenever $X$ is singular, it fails. 
For example, for $n=1$, $L_1 X$ is the tangent space of $X$ and hence
not a locally trivial bundle over $X$.)
The idea of Kontsevich is to give $L_\infty X$
a measure which takes values in the Grothendieck ring $\MM$
of $k$-varieties
which is localized by the class $\LL$  of the affine line.
For each $n \in \Zz$, the family of constructible subsets of $L_n X$ is
 stable under finite union or finite intersection. In other words,
this family is
a Boolean algebra. 
The map
\begin{align*}
  \{ \text{constructible subsets of } L_n X \} & \to \MM  \\
                    A  \hspace{10mm}  &   \mapsto \{ A \} \LL^{-nd}
\end{align*}
is a finite additive measure.
For each $m>n \in \Zz$, because the truncation
morphism $\pi_n^m:L_m X \to L_n X$ is a locally trivial affine space bundle of
relative dimension $(m-n)d$, the pull-back
\[
(\pi_n^m)^{-1}: \{ \text{constructible subsets of } L_n X\}  \to
     \{ \text{constructible subsets of } L_m X\},
\]
 is considered to be 
an extension of the measure into a bigger Boolean algebra.
The \textit{motivic measure} on $L_\infty X$
 is defined to be the limit of these extensions.
Denef and Loeser generalized the motivic measure
to the case where $X$ is singular \cite{germs}.

The integral
 of a function with respect to the motivic measure produces a new 
invariant. In particular, when $X$ is smooth and the function is constant 
equal to 1,
 then the integral, which is the full volume of
$L_\infty X$,  is the class of $X$ in $\MM$.
It
reduces to the Hodge structure of the cohomology of $X$ if
$X$ is complete.
By the transformation rule of Kontsevich, 
for any resolution $Z \to X$,
it equals
the integral of a function on $L_\infty Z$ 
determined by the relative canonical divisor $K_{Z/X}$.
This implies the following theorem of Kontsevich, which we will
generalize:
\begin{thm}
Let $X$ and $X'$ be smooth complete varieties. 
Suppose that there are proper birational morphisms $Z\to X$ and
$Z\to X'$ such that $K_{Z/X}=K_{Z/X'}$.
Then the rational cohomologies of $X$ and $X'$ have the same Hodge structures.
\end{thm}

It is key that 
by the valuative criterion for properness,
almost every $\infty$-jet on $X$
lifts  to a unique  $\infty$-jet on $Z$, 
and hence the map $L_\infty Z \to L_\infty X$ is bijective outside of
measure zero subsets.
Note that  Batyrev first proved the equality of Betti numbers
in the case where $X$ and $X'$ are Calabi-Yau varieties,
with  $p$-adic integration and the Weil conjecture \cite{CYbetti}.

Let $X$ be a variety with Gorenstein canonical singularities.
Denef and Loeser gave $L_ \infty X$ another measure, called the
\textit{motivic Gorenstein measure}, denoted $\muGor _X$ \cite{DL-quotient}.
As in the case with $X$ smooth, $\muGor _X ( L_ \infty X )$ is calculated by
the relative canonical divisor $K_{Z/X}$ for a resolution $Z \to X$.
This implies:

\begin{prop}
Let $X$ and $X'$ be varieties with Gorenstein canonical singularities.
Suppose that there are proper birational morphisms $Z\to X$ and
$Z\to X'$ such that $K_{Z/X}=K_{Z/X'}$.
Then $\muGor _X ( L_ \infty X ) = \muGor _{X'} ( L_ \infty X' )$
\end{prop}

Quotient singularities form one of the mildest classes of 
singularities\footnote{Here the words `quotient singularities'
mean `quotient singularities with respect to the \'etale topology', see 
Definition \ref{quot-sing}.}.
If $X$ is a variety with quotient singularities, then we can give $X$ an
orbifold structure. In the algebro-geometric context, there is 
a smooth Deligne-Mumford stack $\X$ such that $X$ is the coarse moduli
space of $\X$ and the automorphism group of general points of $\X$ is
trivial. 
Although the natural morphism $\X \to X$
is proper and birational,  \textit{not} almost every $\infty$-jet on $X$
lifts  to a  $\CC[[t]]$-point of $\X$ from lack of the \textit{strict}
valuative criterion for properness.
However by twisting the source $\Spec \CC[[t]]$, 
we can lift almost every $\infty$-jet on $X$ to $\X$.
More precisely, 
A \textit{twisted $\infty$-jet} on $\X$ is a representable morphism 
$\D \to \X$ such that $\D$ is a smooth Deligne-Mumford stack with
 the coarse moduli space $\Spec \CC[[t]]$ and $\D$ contains $\Spec \CC((t))$ as
open substack.
(A paper of Abramovich and Vistoli \cite{AV} 
was the inspiration for this notion.
They introduced the notion of twisted stable map.)
For almost every $\infty$-jet $ \gamma : \Spec \CC[[t]] \to X$,
there is a unique twisted $\infty$-jet $\D \to \X$ such that 
the induced morphism  $\Spec \CC[[t]] \to X$
of the coarse moduli spaces 
 is $\gamma$.
If $\overline{\cL _\infty \X}$ is the coarse moduli space of
the twisted $\infty$-jets on $\X$, then 
we define the \textit{motivic measure} $\mu_\X$ on 
 $\overline{\cL_\infty \X}$ in a similar fashion as on $ L_\infty  X$,
though it takes values in the Grothendieck ring of Hodge structures.
We show the following
 close relation between $\muGor_X$ and $\mu_\X$:
\begin{thm}\label{main-in-intro}
The following equality holds:
\[
\muGor_X = \sum_{\Y \subset I(\X)} \LL^{s(\Y)} \mu_\X .
\]
For the precise meaning, see Theorem \ref{main}.
\end{thm}

Chen and Ruan defined the \textit{orbifold cohomology}
 for arbitrary orbifold \cite{CR}.
It originates from string theory on orbifolds \cite{DHVW}.
Let $X$ be a variety with Gorenstein quotient singularities and
$\X$ as above.
The \textit{inertia stack} of $\X$, denoted $I(\X)$,
 is an object in the algebro-geometric realm that corresponds
 to the twisted sector. 
We define the $i$-th orbifold cohomology group of $X$ by:
\[
 \Horb^i (X, \QQ) 
    := \bigoplus_{ \Y \subset I(\X) } H^{i - 2s(\Y)} (\overline \Y , \QQ )
           \otimes \QQ( -s( \Y ) ),
\]
where $\Y$ runs over the connected components of $I(\X)$,
$ \overline{\Y} $ is the coarse moduli space of $\Y$ and
$s(\Y)$ is an integer which is representation-theoretically determined.

\begin{rem}
\begin{enumerate}
\item
The author guesses, though is not sure,
 that our orbifold cohomology is equal to
one by Chen and Ruan in \cite{CR}.
The point he wonders is whether the cohomology groups of $\overline \Y$ are
isomorphic to
 ones of the analytic orbifold (V-manifold) associated to $\Y$.

\item
The orbifold Hodge numbers (that is, the Hodge numbers associated to
$\Horb ^{i}(X , \QQ)$) are 
equal to Batyrev's stringy Hodge numbers \cite{non-archi}.
It follows from Lemma \ref{Goren} and 
Theorem \ref{main-in-intro}.

\end{enumerate}
\end{rem}

Theorem \ref{main-in-intro}
 implies the fact that when $X$ is complete,
 the invariant $ \muGor_X (L_\infty X) $ 
reduces to the alternating sum of the orbifold cohomology groups of $X$.
Hence we obtain the following, conjectured by Ruan \cite{ruan}:

\begin{thm}[=Corollary \ref{cor-of-main}]\label{thm-in-intro}
Let $X$ and $X'$ be complete varieties with Gorenstein quotient 
singularities.
Suppose that there are proper birational morphisms $Z\to X$ and
$Z\to X'$ such that $K_{Z/X}=K_{Z/X'}$.
Then the orbifold cohomologies of $X$ and $X'$ have the same Hodge structures.
\end{thm}

If $X$ and $X'$ are birational minimal models, that is, 
$ K _ X $ and $ K _ X' $ are nef, then for a common resolution $Z$ of
$ X $ and
$ X' $, we have the equality $ K _ {Z/X} = K _ {Z/X'} $ 
(see \cite[Prop.\ 3.51]{Kollar-Mori}).
Hence $X$ and $X'$ have the same orbifold cohomology with Hodge structure.
Note that in the case where $X$ and $X'$ are global quotients, 
Theorem \ref{thm-in-intro} is due to Batyrev \cite{non-archi} and
Denef-Loeser \cite{DL-quotient}.
After writing out the first version of this paper,
I learned by 
an e-mail message from  Ernesto Lupercio 
 that Mainak Poddar and he
     independently proved Theorem \ref{thm-in-intro}.

\subsection*{ Contents }

The paper is organized as follows.
In Section 1, we review motivic measures. 
Section 2 is the central part of the paper.
Here we introduce the notion of a twisted jet and 
examine the space of them. 
Then we 
prove the main result.
In Section 3, we review Deligne-Mumford stacks and 
prove some general results on Deligne-Mumford stacks which we need in
the preceding section.

\subsection*{ Conventions and Notations }

\begin{itemize}
        \item In Section 1 and 2, we work over $\CC$.
        \item For a Deligne-Mumford stack $\X$, we denote by $\overline \X$
        the coarse moduli space of $\X$.
        \item We denote by $\schS$ (resp. $\sch$) the category of schemes
        over a scheme $S$ (resp. over $\CC$).
        \item For a $\CC$-scheme $X$ (or more generally a stack over $\CC$) and
        a $\CC$-algebra $R$,
         we denote by $X \otimes R$ the product $X \times_\CC \Spec R$.
        Then we denote by $X[[t]]$ (resp.\ $X[[t^{1/l}]]$) the scheme
        $X \otimes \CC[[t]]$ (resp.\ $X \otimes \CC[[t^{1/l}]]$).

\end{itemize}

\subsection*{Acknowledgements}
This paper forms a part of my master's thesis.
I am deeply grateful to Yujiro Kawamata for recommending me to learn 
this area and for his direction and encouragement.
I also would like to
thank James McKernan and Hokuto Uehara for reading this paper and
pointing out English mistakes.


\section{Motivic measure; review}

In this section, we would like to review the theory of motivic measures, 
developed by Kontsevich \cite{Orsay}, Batyrev \cite{non-archi}, 
Denef and Loeser \cite{germs}, \cite{DL-quotient}.
I mention \cite{craw} for a nice introduction and
\cite{looi}, \cite{DL-survey} for surveys.

\subsection{Completing Grothendieck rings}

Let us first construct the ring in which motivic measures take values.

A variety means a reduced scheme of finite type over $\CC$.

\begin{defn}
We define the \textit{Grothendieck ring of varieties}, denoted
$K_0(\mathrm{Var})$, to be 
the abelian group generated by the isomorphism classes $\{X\}$ of 
varieties with the relations $\{X\}=\{X\setminus Y\}+\{Y\}$ if $Y$ is
a closed subvariety of $X$. The ring structure is defined by $\{X\}\{Y\}=
\{X\times Y\}$. 
\end{defn}

In the same fashion, we can define the Grothendieck ring of separated
algebraic spaces of finite type. Actually, it is the same thing
as $K_0(\mathrm{Var})$, since every noetherian algebraic space decomposes
into the disjoint union of schemes \cite[Prop.\ 6.6]{Knutson}.

Suppose that $A$ is a constructible subset of a variety $X$, 
that is, $A$ is a disjont union of locally closed subvarieties $A_i\subset X$.
Then  we put $\{A\}:=\sum_i \{A_i\}\in K_0(\mathrm{Var})$, which is 
independent of the choice of stratifications. 
We denote the class of $\A^1$ by $\LL$ and the localization
$K_0(\mathrm{Var})[\LL^{-1}]$ by $\MM$. For $m\in\ZZ$, let $F_m\MM$ be 
the subgroup of $\MM$ generated by the elements $\{X\}\LL^{-i}$ with
$\dim X -i\le -m$. The collection $(F_m\MM)_{m\in\ZZ}$ is a descending 
filtration of $\MM$ with 
\begin{equation}\label{111}
F_m\MM \cdot F_n\MM\subset F_{m+n}\MM.
\end{equation}

\begin{defn}
We define the \textit{complete Grothendieck ring of varieties} by
\[
\Mhat:= \lim_{\longleftarrow} \MM /F_m\MM.
\]
By condition (\ref{111}), it has a natural ring structure.
\end{defn}
Note that it is not known whether
the natural map $\MM \to \Mhat$ is injective or not.

Recall that a \textit{Hodge structure} is a finite dimensional $\QQ$-vector
space $H$ with a bigrading $H\otimes \CC = \oplus _{p,q\in \ZZ}H^{p,q}$
such that $H^{p,q}$ is the complex conjugate of $H^{q,p}$ and each weight
summand $\oplus_{p+q=m}H^{p,q}$ is defined over $\QQ$.
The category $\mathrm{HS}$ of Hodge structures is an abelian category with
tensor product. 

\begin{defn}
We define 
the \textit{Grothendiek ring of Hodge structures},
 denoted $\GHS$, to be the abelian group which 
consists of formal differences $\{H\}-\{H'\}$, where $\{H\}$ and $\{H'\}$
are isomorphism classes of Hodge structures. 
The addition and the multiplication come from $\oplus$ and $\otimes$
respectively.
\end{defn}

A \textit{mixed Hodge structure} is a finite dimensional $\QQ$-vector
space $H$ with increasing filtration $W_\bullet H$, called the \textit{weight
filtration}, such that the associated graded $\Gr_\bullet^W H$ underlies
a Hodge structure having $\Gr_m^W H$ as weight $m$ summand. 
For a mixed Hodge structure $H$, we denote by $\{H\}$ the element
$\{\Gr_\bullet^WH\}$ of $\GHS$.

The cohomology groups $H_c^i(X,\QQ)$ with compact supports of a variety $X$
has a natural mixed Hodge structure. 

\begin{defn}
We define the \textit{Hodge characteristic} $\chi_h(X)$ of $X$ by
\begin{equation*}
\chi_h(X):=\sum_i (-1)^i \{H_c^i(X,\QQ)\}\in \GHS .
\end{equation*}
\end{defn}
Consider the following map:
\[
(\mathrm{Varieties}) \to \GHS,\, X \mapsto \chi_h (X).
\] 
It factors through the map
\[
(\mathrm{Varieties}) \to \MM ,\, X \mapsto \{X\},
\] because the following hold:
\begin{itemize}
	\item $\chi_h(X\times Y)=\chi_h(X)\chi_h(Y)$,
        \item $\chi_h(X)=\chi_h(X\setminus Y) + \chi_h(Y)$ if 
                     $Y\subset X$ is closed,
        \item the Hodge characteristic  of the affine line,
          $\chi_h(\A^1)= \{H_c^2(\A^1,\QQ)\}$, is invertible. 
\end{itemize}
We also denote by $\chi_h$ the induced homomorphism $\MM \to \GHS$.

For $m\in\ZZ$, let $F_m\GHS$ be the subgroup generated by the elements
$\{H\}$ such that the maximum weight of $H$ is less than or equal
to $ -m $.

\begin{defn}
We define the \textit{complete Grothendieck ring of Hodge structures}, 
denoted $\hat{\GHS}$, as follows:
\[
\hat{\GHS}:= 
\lim_{\longleftarrow} \GHS / F_m\GHS.
\]
\end{defn}
We can see that the natural map $\GHS\to\hat{\GHS}$ is injective. 
Because the maximal weight of $ H^i ( X ,\QQ )$ does not exceed
 $ 2 \dim X $, $\chi_h$ extends to 
$
\chi_h :\Mhat \to \hat{\GHS}.
$


\subsection{Jets on schemes}
For convenience's sake, we denote by $(t^\infty)$ the ideal $(0)$
of the power series ring $\CC[[t]]$.
For $n\in \Zz \cup \{\infty\}$,
we denote by $D_n$ the affine scheme $\Spec \CC[[t]]/(t^{n+1})$.

\begin{defn}
Let $X$ be a scheme. 
For $n \in \Zz \cup \{ \infty \}$, 
we define  the \textit{scheme of $n$-jets}\footnote{
In \cite{germs}, an $\infty$-jet is called an \textit{arc}. 
As it is more convenient, I prefer my terminology.}
of $X$, denoted $L_n X$,
to be the scheme representing the functor
\begin{align*}
\sch &\to (\mathrm{Sets}) \\
U &\mapsto \Hom_{\sch}(U\times D_n, X).
\end{align*}
\end{defn}
Greenberg \cite{greenberg1} proved 
the representability of the functor for $n < \infty$. 
For $m,n \in \Zz$ with $m<n$, 
a canonical closed immersion $D_m \to D_n$ induces a canonical
projection $L_nX \to L_m X$. 
Since all these projections are affine morphisms,
the projective limit 
$L_\infty X = \displaystyle \lim_{\longleftarrow} L_nX$ exists in
the category of schemes.

If $X$ is of finite type, 
then for $n < \infty$, so is $L_n X$.
If $X$ is smooth and of pure dimension $d$, then,
 for each $n \in \Zz$, the natural projection
$L_{n+1}X\to L_nX$ is a Zariski locally trivial $\A^d$-bundle.
If $f:Y\to X$ is a morphism of schemes, then for each $n$,
there is a canonical morphism $f_n:L_n Y \to L_n X$.


\subsection{Motivic measure}\label{motivicmeasure}

Let $X$ be a scheme of pure dimension $d$.
By abuse of notation, we denote the set of points of $L_\infty X$
also by $L_\infty X$. Let $\pi_n : L_\infty X \to L_n X$ be 
the canonical projection.

\begin{defn}\label{def-stable}
A subset $A$ of $L_\infty X$ is \textit{stable at level $n$} if we have:
\begin{enumerate}
	\item $\pi_n(A)$ is a constructible subset in $L_n X$,
        \item $A=\pi_n^{-1}\pi_n(A)$,
        \item for any $m\ge n$, the projection
$
\pi_{m+1}(A) \to \pi_m(A)
$
is a piecewise trivial $\A^d$-bundle.
\end{enumerate}
(A morphism $f:Y\to X$ of schemes is called a \textit{piecewise trivial
$\A^d$-bundle} if there is a stratification $X=\coprod X_i$ such that 
$f|_{f^{-1}(X_i)}:f^{-1}(X_i) \to X_i$ is isomorphic to 
$X_i\times \A^d \to X_i$ for each i.)
A subset $A$ of $L_\infty X$ is \textit{stable} if it is 
stable at level $n$ for some $n\in \Zz$.
\end{defn}

The stable subsets of $L_\infty X$ constitute a Boolean algebra. 
If $A \subset L_\infty X$ is a stable subset, then
$\{\pi_m(A)\}\LL^{-md}\in\Mhat $ is constant for $m \gg 0$.
We denote it by $\mu_X(A)$ \footnote{
This differs from the definition in \cite{craw}, \cite{non-archi} and 
\cite{germs} by a factor $\LL^d$.
}.
The map
\begin{align*}
\mu_X:\{\text{stable subsets of }L_\infty X \}&\to \Mhat \\
                      A\hspace{10mm}&\mapsto \mu_X(A) 
\end{align*}
is a finite additive measure. Let us extend this measure to a bigger
family of subsets of $L_\infty X$.

\begin{defn}
A subset $A \subset L_\infty X$ is called \textit{measurable} if, 
for every $m\in \ZZ$,
there are stable subsets $A_m \subset L_\infty X$ and 
$C_i\subset L_\infty X$, 
$i\in \ZZ_{>0}$
such that the symmetric difference $(A \cup A_m) \setminus (A \cap A_m)$
 is contained in 
$\cup_i C_i$ and we have $\mu_X(C_i) \in F_m \MM$ for all $ i$, and 
$\displaystyle \lim_{i\to \infty} \mu_X(C_i) = 0$ in $\Mhat$.
\end{defn}

The measurable subsets of $L_\infty X$ also constitute a Boolean algebra.
Suppose that $A\subset L_\infty X$ is measurable and $A_m\subset L_\infty X$, 
$m\in\ZZ$ are
stable subsets  as in the definition above.  
We put
$ \mu_X(A):=\displaystyle \lim_{m\to\infty}\mu_X(A_m)$.
It is independent of the choice of $A_m$, see \cite[Prop.\ 2.2]{looi}, 
\cite[Th. A.6]{DL-quotient}.
The map 
\begin{align*}
\mu_X:\{\text{measurable subsets of }L_\infty X \}&\to \Mhat \\
                      A\hspace{10mm}&\mapsto \mu_X(A)
\end{align*}
is a finite additive measure.

\begin{defn}
We call this  the \textit{motivic measure} on $L_\infty X$.
\end{defn}

\begin{defn}
Let $A\subset L_\infty X$ be a measurable subset and
$\nu:A \to \ZZ \cup \{ \infty \}$ a function.
We say that $\nu$ is a \textit{measurable
function} if 
the fibers 
are measurable and $\mu_X(\nu^{-1}(\infty)) = 0$. 
For a measurable function $\nu$, 
we formally define the
\textit{motivic integration} of $\LL^\nu$ by
\[
\int_A \LL^\nu d\mu_X := \sum_{n\in\ZZ}\mu_X(\nu^{-1}(n))\LL^n.
\]
We say that $\LL^\nu$ is \textit{integrable} if the series above converges
in $\Mhat$.
\end{defn}

\begin{expl}
Let  $\I$ be an ideal sheaf on $X$. A point $\gamma \in L_\infty X$ corresponds to a
morphism $\Spec \kappa \to L_\infty X$ for the residue field $\kappa$ of 
$\gamma$
 and hence to a morphism 
$\gamma' : \Spec \kappa[[t]] \to X$. 
The function
\begin{align*}
\ord \, \I:L_\infty X & \to \Zz\cup \{\infty\} \\
\gamma \; & \mapsto  n  \;  \text{ if }(\gamma')^{-1}\I=(t^n) 
\end{align*}
is a measurable function by the following lemma.
\end{expl}

\begin{lem}\cite[Prop.\ 3.1]{looi}, \cite[Lem.\ 4.4]{germs}.
For a subvariety $Y \subset X$ of positive codimension, 
the subset $L_\infty Y \subset L_\infty X$ is of measure zero.
\end{lem}

\begin{expl}\label{formula1}\cite[Th.\ 3.6]{non-archi}, \cite[Th.\ 2.15]{craw}.
Let $X$ be a smooth variety of dimension $d$ 
and $E = \sum_{i=1}^{r} d_i E_i$ an effective divisor
on $X$ with simple normal crossing support.
For a subset $J \subset \{ 1, \dots , r \}$, we define
\[ 
E _J ^\circ := \bigcap _{i \in J} E _i  \bigg \backslash 
\bigcup _{ i \in \{  1, \dots , r  \} \setminus J}  E _i .
\]
If $\I _E$ is the ideal sheaf associated to $E$, then we have the
following formula:
\[
   \int _{ L_\infty X  } \LL ^{ - \ord \, \I_E } d \mu _X
 = \sum _{ J \subset \{ 1, \dots , r   \}  } \{ E _J ^\circ \} 
     \prod _{ i \in J  } \frac{ \LL - 1 }{ \LL ^{ d_i + 1 } - 1 } \ .
\]
\end{expl}

\subsection{The transformation rule}

Let $X$ and $Y$ be  varieties of pure dimension $d$ and 
$f:Y[[t]] \to X[[t]]$ a morphism over $D_\infty$.
We define the morphisms $f_n:L_n Y \to L_n X $ as follows.
For a scheme $U$, a $U$-point $\gamma$ of $L_n Y$ is a morphism
$\gamma':D_n \times U \to Y$.  If  $\psi$  is the composition of
natural morphisms
$D_n \times U \xrightarrow{\pr} D_n \to D_\infty$, then
we define $f_n(\gamma)$ to be the $U$-point of $L_n X$ 
corresponding to the composition
\[
  D_n \times U \xrightarrow{(\gamma',\psi)} 
 Y \times D_\infty = Y[[t]]  
                 \xrightarrow{f} X[[t]] \to X.
\]

Assume that $Y$ is smooth.
We define the \textit{jacobian ideal sheaf} $\J_f$ of $f$ to be the $0$-th
Fitting ideal sheaf of $\Omega_{Y[[t]]/X[[t]]}$ 
(see \cite[V.1.3]{Eisen-Harris} or \cite[20.2]{Eisen}).
This is the ideal sheaf such that 
$\J_f \Omega_{Y[[t]]/D_\infty}^d =
f^* (\Omega^d_{X[[t]]/D_\infty}/(tors))$
as subsheaves of $\Omega_{Y[[t]]/D_\infty}^d$, where $(tors)$ is the 
torsion.
The following, called the \textit{transformation rule},
 is the most basic theorem in the theory.

\begin{thm}
\cite[Th.\ 1.16]{DL-quotient}, \cite[Th.\ 3.2]{looi}\label{changevariables}.
Let $A$ be a measurable set in $L_\infty Y$ and 
$\nu:f_\infty(A)\to \ZZ \cup \{\infty\}$ a measurable function.
Suppose that $f_{\infty}|_A$ is injective.
Then we have the following equality:
\[
\int_{f_\infty (A)} \LL^\nu  d\mu_X = 
  \int_A \LL^{\nu \circ f_\infty - \ord \, \J_f} d\mu_Y.
\]
\end{thm}

We will generalize this later (Theorem \ref{thm-GeneralTrans}).


\subsection{The motivic Gorenstein measure}\label{motivic-Gor}

Let $X$ be a variety with 1-Gorenstein and canonical singularities,
 that is, the canonical sheaf
$\omega_X$ is invertible and 
all discrepancies are $\ge 0$, (see \cite[\S 0-2]{KMM}).
Then there exists a natural morphism 
$\Omega_X^d \to \omega_X$. The kernel of this morphism is
the torsion. 
We define an ideal sheaf $\I_X$ on $X$ by the equation:
\[
\I_X \omega_X = \Img (\Omega_X^d \to \omega_X) .
\] 
Then $\LL^{\ord \, \I_X}$ is integrable by Example \ref{formula1} and 
Lemma \ref{Goren}.

\begin{defn}
We define the \textit{motivic Gorenstein measure} $\muGor_X$ on $L_\infty X$ 
as follows:
\begin{align*}
\muGor_X:\{\text{measurable subsets of }L_\infty X \}&\to \MM \\
                      A\hspace{10mm}&\mapsto \int_A \LL^{\ord \, \I_X} d\mu_X .
\end{align*}
\end{defn}

\begin{lem}\label{Goren}
Let $X$ and $X'$ be complete varieties with 1-Gorenstein
canonical singularities.
\begin{enumerate}
   \item Let $A$ be a measurable subset of $L_\infty X$ and $f:Z \to X$ be
       a resolution. Then 
       \[
            \int_A \LL^{\ord \, \I_X} d\mu_X = 
             \int_{f_\infty^{-1}(A)} \LL^{-\ord \, \cK} d \mu_Z , 
        \]
      where $\cK$ is the ideal sheaf associated to $K_{Z/X}$.
   \item Suppose that there exist proper birational morphisms 
        $f:Z\to X$ and
        $g:Z\to X'$ with $K_{Z/X}=K_{Z/X'}$. Then we have 
       $\muGor_X(L_\infty X)= \muGor_{X'}(L_\infty X')$.
\end{enumerate}
\end{lem}

\begin{proof}
(1). By Theorem \ref{changevariables}, we have 
\begin{align*}
\int_A \LL^{\ord  \I_X} d\mu_X  
      = \int_{f_\infty^{-1}(A)} \LL^{\ord  \I_X \circ f_\infty
          - \ord  \J_f} \mu_Z .
\end{align*}
We have to show 
\begin{equation}\label{777}
\ord {\I_X}\circ f_\infty - \ord  \J_f = - \ord  \cK.
\end{equation}

Pulling back $\I_X \omega_X \cong \Omega_X^d/
(tors)$, we have $(f^{-1}\I_X) (f^* \omega_X)
 \cong \J_f \Omega_Z^d$.
On the other hand, we have 
$f^* \omega_X \cong \cK \omega_Z \cong \cK \Omega_Z^d$.
Hence $(f^{-1}\I_X) \cdot  \cK = \J_f$. This shows (\ref{777}).

(2) is a direct consequence of (1).
\end{proof}

\begin{rem}
The invariant $\muGor_X(L_\infty X)$ has been
 already introduced in \cite{Orsay}
by using a resolution of singularities. 
Denef and Loeser constructed the invariant more directly with the motivic
Gorenstein measure \cite{DL-quotient}.
\end{rem}

Suppose $X$ is complete.
The question is whether $\chi_h \circ \muGor_X(L_\infty X)$ is 
the alternating sum of \textit{a kind of cohomology groups}, as the case 
where $X$ is smooth. It is known that when $X$ is a global quotient, the
answer is Yes \cite{non-archi}, \cite{DL-quotient}. 
Our result, Theorem \ref{main}, says that when $X$ has only quotient
singularities, the answer is Yes.


\section{Twisted jets}

In this section, we need manage the theory of Deligne-Mumford stacks.
See the next section for the generalities about Deligne-Mumford stack.

\subsection{Non-twisted jets on stacks}

The following is a direct generalization of the notion of
jet on schemes:

\begin{defn}
Let $\X$ be a Deligne-Mumford stack. 
For $n\in\Zz\cup \{\infty\}$,
we define the \textit{stack of non-twisted $n$-jets} of $\X$, denoted 
$L_n\X$, as follows.
An object of $L_n\X$ over $U\in \sch$ is an object of $\X$ over $U\times D_n$.
 For a morphism $\varphi:V\to U$ in $\sch$,
a morphism in $L_n \X$ over $\varphi$ is a morphism in $\X$ over
$\varphi \times \id_{D_n}$.
\end{defn}

\begin{lem}
For every $n\in \Zz \cup \{\infty\}$, 
$L_n \X$ is a stack. 
\end{lem}

\begin{proof}
It is clear that they satisfy the axioms of category fibered in groupoids.
If $(U_i\to U)_i$ is an \'etale covering in $\sch$, then 
so is $(U_i\times D_n \to U\times D_n)_i$.
As a result, they satisfy the axioms of stack, too.
\end{proof}

Let $f:Y \to X$ be a morphism of schemes and $Z\subset Y$ a closed 
subscheme with an ideal sheaf $\mathfrak a$.
We say that $f$ is $Z$-\textit{\'etale}
 if for any ring $A$ and any nilpotent ideal $J \subset A$, and for
any commutative diagram of solid arrows
\[
\xymatrix{
  \Spec A/J \ar[r]^{\varphi} \ar[d] & Y \ar[d] \\
  \Spec A \ar[r] \ar@{-->}[ur] & X 
}
\]
such that $\varphi^{-1} \mathfrak a $ is nilpotent, 
there is a unique broken arrow
which makes the whole diagram commutative.
$Z$-\'etaleness is defined for a representable morphism of stacks also in
the evident fashion.

\begin{lem}\label{jet-atlas}
\begin{enumerate}
	\item Let $M$ be a scheme and $N \subset M$ a closed subscheme.
        We denote by $(L_nM)_N$ the subscheme of $L_n M$ parametrizing
       the jets with the base point in $N$.
          Let $p:M\to \X$ be an 
           $N$-\'etale morphism.
             Then, for every $n\in \Zz \cup \{\infty\}$, 
            we have a natural isomorphism:
           \[
              (L_nM)_N \cong L_n\X \times_\X N.
            \] 
            \item For every $n\in \Zz \cup \{\infty\}$,
            $L_n \X$ is a Deligne-Mumford stack.
\end{enumerate}
\end{lem}

\begin{proof}
(1). Let us first show $L_n\X \times_\X N$ is an algebraic space.
Let $\pi:L_n\X \to \X$ be the canonical projection. An object of 
$L_n\X\times_\X N$ is a triple $(\gamma, f, \alpha)$ 
where $\gamma : U \times D_n \to \X $,
$f:U \to N$, and $\alpha:\pi(\gamma)\to p(f)$ is a morphism
in $\X$ over $U$.
 By definition, an automorphism of $(\gamma, f, \alpha)$ is an 
automorphism $\theta$ of $\gamma$ such that $\alpha \circ \pi(\theta)
=p(\id_f)\circ\alpha$. Hence $\theta$ must be the identity. 
We have thus proved that the automorphism of every object
of $L_n\X \times_\X N$ is trivial and hence that $L_n\X \times_\X N$ is 
an algebraic space (see \cite[Cor.\ 8.1.1]{LMB}).

The diagram of solid arrows
\[ \xymatrix{
U \ar[r]^f \ar@{^{(}->}[d] & N \ar@{^{(}->}[r] &  M  \ar[d]^p \\
U\times D_n \ar[rr]_\gamma \ar@{-->}[urr]^\tau & & {\X}
} \]
is commutative. 
We can see that there is a unique broken arrow $\tau$ in
the diagram. If $n < \infty$, since $p$ is $N$-\'etale, 
this is trivial from the definition.
If $n = \infty$, since $D_ {\infty}$ is the direct limit
of $  D_n $, $0 \le  n < \infty$, this follows from
 the case $n < \infty$. 
Sending $(\gamma, f, \alpha)$ to $\tau$ defines a morphism
$ \Phi: L_n\X \times_\X M \to L_nM$. 

The inverse of the morphism $\Phi$  is given by
$(p_n,\pi_M): (L_nM)_N \to L_n\X \times_\X N $ where $p_n:L_nM \to L_n\X$
and $\pi_M:(L_n M)_N \to N=(L_0 M)_N$ are the natural morphisms.
We have thus proved (1).

(2). Now suppose that $p$ is \'etale and surjective.
Consider the following cartesian diagram:
\[\xymatrix{ \ar @{} [dr] |{\square}
L_nM \ar[r] \ar[d]_{\pi_M} & L_n \X \ar[d]^{\pi} \\
 M  \ar[r]_p  &  \X  .
}\]
Because $\pi_M$ is representable and $p$ is \'etale and surjective,
$\pi$ is representable (see \cite[Lem.\ 4.3.3]{LMB}).
This completes the proof
(see \cite[Prop.\ 4.5]{LMB}).
\end{proof}


\subsection{Twisted jets}
For a positive integer $l$,
we put $\zeta_l := \exp (2\pi \sqrt{-1} / l)$.
 Let $\mu_l= \langle \zeta_l \rangle$ be
the group of the $l$-th roots of $1$. 
$\mu_l$ acts on $D_n$ by $\zeta_l: t \mapsto \zeta_l t$.
We denote by $\D^l_n$ the quotient stack 
$[D_m/\mu_l]$ with $m=nl$. 
The stack $\D^l_n$ has the canonical atlas $D_m \to \D^l_n$ and the closed
 point $\Spec \CC \to \D^l_n$. 
We fix a morphism $ \D^l_n \to D_n $ such that $ D_n $ is the coarse
moduli space of $ \D^l_n $ for this morphism, and such that the composition 
\[
     D_m \to \D^l_n \to D_n
\]
is given by the ring homomorphism $ \CC[[t]]/(t^{n+1}) \to  \CC[[t]]/(t^{m+1})$,
$ t \mapsto t^l $. 
\begin{defn}
Let $\X$ be a Deligne-Mumford stack. 
A \textit{twisted $n$-jet of order} $l$
 on $\X$ is 
a representable
morphism $\D^l_n \otimes \Omega \to \X$
for an algebraically closed field $\Omega \supset \CC$.
\end{defn}

For a Deligne-Mumford stack $\X$,
the \textit{inertia stack} of $\X$, denoted $I(\X)$, 
is the stack parametrizing the pairs $(\xi,\alpha)$
such that $\xi\in \ob\X$ and $\alpha \in \Aut(\xi)$.
For details on inertia stack, see Subsection \ref{subsec-inertia}.
There is a natural forgetting morphism $I(\X)\to \X$.
For $l\in \ZZ_{>0}$, let $I^l(\X) \subset I(\X)$ denote the substack
parametrizing the pairs $(\xi, \alpha)$ with $\ord( \alpha) = l$.

Let $\gamma:\D^l_n \otimes \Omega \to \X$ be a twisted $n$-jet 
of order $l$ on $\X$.
The canonical morphism
\[ 
\tilde \gamma:D_m \otimes \Omega \to \D^l_n \otimes \Omega \to \X
\]
is considered to be an $\Omega$-point of $L_m \X$ and 
the canonical morphism 
\[ 
\overline \gamma:\Spec \Omega \to D_m \otimes \Omega \to \D^l_n \otimes \Omega 
\to \X
\]
to be an $\Omega$-point of $\X$. 
Since the automorphism group of the closed point of $\D^l_n \otimes \Omega$
is identified with $\mu_l$, $\gamma$ induces an injection 
$\mu_l \to \Aut (\overline \gamma)$. 
If $b\in \Aut (\overline \gamma)$ is the image of $\zeta_l$, then the pair
$(\overline \gamma, b)$ is regarded as an $\Omega$-point of $I^l(\X)$
and the triple $(\tilde \gamma,(\overline \gamma, b), \id_{\overline \gamma})$
 to be an $\Omega$-point of
$L_m \X \times_\X I^l(\X)$.
We define a map $\Psi$ by
\begin{align*} 
\Psi :\{ \text{twisted n-jets  of order $l$ on }\X \} & \to 
|L_m \X \times_\X I^l(\X)|  \\
 \gamma \hspace{10mm} & \mapsto 
       (\tilde \gamma,(\overline \gamma, b), \id_{\overline \gamma}) .
\end{align*}

\begin{lem}\label{closed}
The subset $\Img(\Psi)\subset|L_m \X \times_\X I^l(\X)| $ is closed 
for the Zariski topology.
\end{lem}

\begin{proof}
Fix an atlas $p:M \to \X$ with $M$ separated.

We will first characterize the points in $\Img(\Psi)$.
On account of the arguments on groupoid spaces in Subsection \ref{subsec-DM}, 
we can see that the following are equivalent:
\begin{enumerate}
	\item to give a commutative diagram
        \[\xymatrix{ 
           D_m\otimes \Omega \ar[d] \ar[r] ^ \eta & M \ar[d]^p  \\
             \D^l_n \otimes \Omega \ar[r]_{\gamma}  &   \X
        }\]
        such that $\gamma$ is a twisted $n$-jet of order $l$,
        \item to give a morphism of groupoid spaces
        \begin{equation*}\label{diagram}\xymatrix{
         {(D_m\otimes \Omega)\times \mu_l} \ar[r] ^ { \delta' }
        \ar@<-1ex>[d]_\pr \ar@<1ex>[d]^{\mu_l-\text{action}} & 
       M\times_\X M \ar@<-1ex>[d]_{\pr_1} \ar@<1ex>[d]^{\pr_2}\\
        D_m\otimes \Omega \ar[r] _ \eta & M
       }\end{equation*}
        such that the composition
        \[
           \Spec \Omega \hookrightarrow (D_m\otimes \Omega) \times 
           \{ \zeta _ l \} \hookrightarrow (D_m\otimes \Omega)\times \mu_l
           \xrightarrow { \delta' } M\times_\X M
        \]
         corresponds to an automorphism of order $l$ of the following
        $ \Omega $-point of $\X$:
        \[
            \Spec \Omega \hookrightarrow D_m\otimes \Omega 
                 \xrightarrow { \eta } M \xrightarrow { p } \X .
        \]
        \item  to give a morphism $\delta:D_m \otimes \Omega \to M\times_\X M$
          such that $\pr_1 \circ \delta  \circ  \zeta_l=\pr_2 \circ \delta$
          and  the composition
        \[
           \Spec \Omega \hookrightarrow D_m\otimes \Omega
           \xrightarrow { \delta } M\times_\X M
        \]
         corresponds to an automorphism of order $l$ of the following
        $ \Omega $-point of $\X$:
        \[
            \Spec \Omega \hookrightarrow D_m\otimes \Omega 
                 \xrightarrow { \delta } M \times _ \X M
              \xrightarrow { \pr _ 1 } M \xrightarrow { p } \X .
        \]
\end{enumerate}
Any point of $ | L_m \X \times_\X I^l(\X) | $ is represented by 
the triple $ (\psi, (\overline \psi, b ) , \id _{ \overline \psi }) $
such that $ \psi $ is an $\Omega$-point of $L_m \X$ with 
an algebraically closed field $\Omega$, 
$\overline \psi$ is an $\Omega$-point of 
$\X$ corresponding to the composition
\[
 \Spec \Omega \hookrightarrow D_m \otimes \Omega \xrightarrow{ \psi} \X
\]
 and $b$ is an automorphism of $\overline \psi$.
Then the equivalence above implies: 
\begin{list}{$\clubsuit$}{}
	\item $ (\psi, (\overline \psi, b ) , \id _{ \overline \psi }) $
 is in $\Img(\Psi)$ iff for a lift 
$\eta :D_m \otimes \Omega \to M$ of $\psi$, there exists a morphism
$\delta:D_m \otimes \Omega \to M\times_\X M$ such that 
$\pr_1 \circ \delta = \eta$ and 
$\pr_2 \circ \delta = \eta \circ  \zeta_l$
and the composition $\Spec \Omega \to D_m \otimes \Omega \xrightarrow{\delta}
 M\times_\X M$ 
corresponds to $b$.
\end{list}

Let $\xi$ and $\sigma$ be points of $L_m \X \times_\X I^l(\X)$.
Suppose that $\sigma$ is in $\Img (\Psi)$ and $\xi$ is a specialization of 
$\sigma$.
It suffices to show that $\xi$ is in $\Img (\Psi)$.
By \cite[Prop.\ 7.2.1]{LMB}, there is a complete discrete valuation 
ring $R$ with algebraically closed
 residue field $\kappa$ and quotient field $K$ such that 
there is a commutative diagram as follows:
\[\xymatrix{
{\Spec K} \ar[d] \ar[dr]^\sigma      &            \\
{\Spec R} \ar[r]^<<<<{\theta}            &     L_m \X \times_\X I^l(\X)   \\
{\Spec \kappa} \ar[u] \ar[ur]_\xi . & 
}\]
(Here by abuse of notation, the arrows $\sigma$ and $\xi$ in the diagram 
are representatives of the points $\sigma$ and $\xi$ respectively.)
If $ \theta $ corresponds to a triple
$(\lambda_R, (\overline {\lambda_R}, b_R ) , \id _{ \overline {\lambda_R} })$,
then, the pull-backs 
$(\lambda_K, (\overline {\lambda_K}, b_K ) , \id _{ \overline {\lambda_K} })$
and
$(\lambda_\kappa, (\overline {\lambda_\kappa}, b_\kappa ) ,
 \id _{ \overline {\lambda_\kappa} })$ correspond to
$\sigma$ and $\xi$ respectively.
By extending $R$, we can assume that 
$\overline {\lambda_R}: \Spec R \to \X$ lifts to 
 $\nu_R : \Spec R \to M$.
Since $p$ is \'etale,
$\lambda_R : D_m \otimes R \to \X $
uniquely lifts to $\tilde \nu :D_m \otimes R \to M$ such that the diagram
\[\xymatrix{
            \Spec R  \ar@{^{(}->}[d] \ar[r]^\nu    &     M  \ar[d]  \\
          D_m \otimes R \ar[r]_(.6){\lambda_R} \ar[ur]^{\tilde \nu} &   \X
}\]
is commutative.
 Let $\overline K$ be the algebraic closure of $K$,
let $\eta$ be the composition 
$D_m \otimes \overline K \to D_m \otimes R \xrightarrow{\tilde \nu} X$ and
let  $b_R': \Spec R \to M \times_\X M$ be the lift of $\nu_R$ which
corresponds to $b_R$.
From $\clubsuit$ and the assumption, 
there is a morphism 
$\delta:D_m\otimes \overline K \to M\times_\X M$ such that 
$\pr_1 \circ \delta = \eta$ and 
$\pr_2 \circ \delta = \eta \circ \zeta_l $, and 
the composition 
\[
 \Spec \overline K \hookrightarrow D_m\otimes \overline K \xrightarrow{\delta}
 M\times_\X M
\] 
equals the composition
\[
 \Spec \overline K \to \Spec R \xrightarrow{b'_R}  M\times_\X M .
\]
We can replace $\overline K$ with a finite extension $K'$ of $K$. 
Moreover, replacing $R$ with its normalization in $K'$, we can assume that
$K'=K$ and that $b'_R$ and $\delta$ 
induce the same morphism $\Spec K \to M\times_\X M$.
\[\xymatrix{
{\Spec K} \ar[r] \ar@{^{(}->}[d] & {\Spec R} \ar[r]^{b'_R} \ar@{^{(}->}[d] 
                           & M\times_\X M \ar@<-1ex>[d] \ar@<1ex>[d] \\
D_m\otimes K \ar[r] \ar[urr]^(.3){\delta} & D_m\otimes R \ar[r]_{\tilde \nu} & M
}\]

Consider the unique morphism $\tau :D_m\otimes R \to M\times_\X M$ 
such that $ \pr_1 \circ \tau  = \tilde \nu$.
Then the two morphism  
$\pr_2\circ \tau $ and $\tau\circ\zeta_l$ is the same morphism because of
the separatedness of $M$. Then the composition $D_m \otimes \kappa \to
D_m \otimes R \xrightarrow{\tau}
 M\times_\X M$ satisfies the condition in $\clubsuit$.
Hence $\xi \in \Img (\Psi)$. The proof is now complete.
\end{proof}

\begin{defn}
We define the \textit{stack of twisted $n$-jets of order $l$} on $\X$,
 denoted 
$\cL_n^l \X$, to be the reduced closed substack of $L_m \X \times_\X I^l(\X)$
with support $\Img(\Psi)$. 
We define the \textit{stack of twisted $n$-jets} on $\X$,
denoted $\cL_n\X$,
to be the disjoint sum  
$\coprod_{l\ge 0} \cL_n^l\X $.
In  particular,  $\cL_0 \X$
  is the inertia stack $I(\X)$.
\end{defn}

If we set 
\[
l_0 := \max \{l \ | \ l= \ord (\alpha), 
\text{ for some } \xi \in \ob \X \text{ and for some } 
\alpha \in \Aut (\xi) \},
\]
then for any $l > l_0$,  $\cL_n^l\X = \emptyset$.
So the disjoint sum above is indeed a finite sum.

\subsection{The formal neighborhood of $I(\X)$ and its canonical automorphism}

Let $\X$ be a smooth 
Deligne-Mumford stack, $x:\Spec k \to \X$ its closed point.
Then the tangent space $T_x \X$ is defined to be $T_{v} M$ for an atlas
$M \to \X$ and a lift $ v:\Spec k \to \X $ of $x$, uniquely determined up to
unique isomorphism. 
Then $\Aut (x)$ naturally acts on $T_x \X$. 
We now globalize it.

Let $\Y$ be a connected component of the inertia stack
$I(\X)$, $F: \Y \to \X$ the 
forgetting map. 
Let $\X_0$ be the image of $\Y$ by $F$. 
The completion $\hat{\cO_\X}$ of $\cO_\X$ along $\X_0$ is considered to be 
an $\cO_{\X_0}$-algebra.
Then we define a coherent sheaf $\cA$ on $\Y$ to be the pull-back
of $\hat{\cO_\X}$ by $F:\Y \to \X_0$.

\begin{defn}
We define $\cN := \cSpec \cA$ and call it the \textit{formal neighborhood}
 of $\Y$.
\end{defn}

If we set $\hat{\X} := \cSpec \hat{\cO_\X}$ where we consider
$\hat{\cO_\X}$ to be an $\cO_{\X}$-algebra, there is a natural
morphism $\hat{\X} \to \X$, which is $\X_0$-\'etale. 
Since $\Y \to \X_0$ is unramified and flat, it is \'etale.
Hence
 the natural morphism
$\cN \to \hat{\X}$ is also \'etale and the composition 
$\cN \to \hat{\X} \to \X$ is $\Y$-\'etale.

Let $U$, $V$ be varieties.
Let $\xi: U \to \X$ be an \'etale morphism, $\nu : V \to U$ 
a morphism and $\alpha$ an automorphism of $ \xi \circ \nu $.
Suppose that $ (\xi \circ \nu, \alpha): V \to \Y $ is \'etale.
\[
\xymatrix{
 V \ar[d]_{(\xi \circ \nu ,\alpha)} \ar[r]^\nu & U \ar[d]^\xi \\
 \Y \ar[r] & \X
}
\]
Then we obtain a commutative diagram
\[\xymatrix{
 V \ar[r]^\nu \ar[d]_\nu &  U \ar[d]^\xi \ar@{=>}[dl]_\alpha \\
 U  \ar[r]_\xi &  \X ,
}\]
where `$\Rightarrow$' denotes a 2-morphism. 
Let $ \tilde \alpha : V \to U \times _\X U $ be the corresponding morphism.
\[\xymatrix{
      & U\times_\X U  \ar@<-1ex>[d]_ {\pr _1} \ar@<1ex>[d]^{\pr_2} \\
      V \ar[r]_\nu \ar[ur]^{\tilde \alpha}     & U .
}\]
If $\hat \cO_U$ is the completion of $\cO_U$ along $\nu( V)$, 
then $\cA_V = \nu^* \hat \cO_U$.
 We have a canonical automorphism of $\cA_V$
\[
  \cA_V \xrightarrow{\pr_2^*} \tilde \alpha^* \hat \cO_{U\times_\X U} 
      \xrightarrow{(\pr_1^*)^{-1}} \cA_V ,
\]
and hence a canonical automorphism of $\cA$ and $\cN$.
Now this automorphism of $\cN$ is considered to be a globalization of
the action on $T_x \X$ mentioned above.

\subsection{Shift number and orbifold cohomology}

Suppose that $\Y$ is contained in $I^l(\X)$ for an integer 
$l \ge 1$.
Let $(x,\alpha)$ be a closed point of $\Y$ where $x$ is a closed point of 
$\X$ and $\alpha \in \Aut (x)$. 
Then $\alpha$ acts on the tangent space $T_x \X$.
For a suitable basis, this automorphism is given
by a diagonal matrix 
\[
\diag (\zeta_l^{a_1}, \dots , \zeta_l^{a_d})
\]
 with $1 \le a_j \le l$ and $d= \dim \X $.

\begin{defn}
We define the \textit{shift number} of $\Y$ by 
\[
    s(\Y) := \dim \X - \frac{1}{l} \sum _{j=1}^d a_j
              = \frac{1}{l} \sum _{j=1}^d (l-a_j).
\]
\end{defn}

This is determined by the rank of the eigenbundles of 
$\cN$ for the canonical action. Hence it depends only on $\Y$.

Suppose that the coarse moduli space $X=\overline \X$
is a variety with Gorenstein quotient singularities and 
$\X$ has no reflections.
Then the matrix $\diag (\zeta_l^{a_1}, \dots , \zeta_l^{a_d})$
is in $\SL_d(\CC)$ (see \cite{K-Watanabe}). Hence $s(\Y)$ is an integer.

Now, let us define the orbifold cohomology.
\begin{defn}\label{orb-cohomology}
Assume $X$ is complete.
Then we define the \textit{$i$-th orbifold cohomology group}
 along with Hodge structure
as follows:
\begin{equation*}
\Horb^i(X,\QQ):=\bigoplus _\Y H^{i-2s(\Y)}(\overline{\Y},\QQ)
       \otimes \QQ(-s(\Y)),
\end{equation*}
where $\Y$ runs over the connected components of $I(\X)$
and $\QQ(-s(\Y))$ is a Tate twist $\QQ(-1)^{\otimes s(\Y)}$.
\end{defn}

Since the natural morphism $\overline \Y \to X$ is quasi-finite,
$\overline \Y$ is a scheme (see \cite[Th.\ A.2]{LMB}).
Because of this and Corollary \ref{cor-inertia-smooth},
 $\overline \Y$ is a complete variety with 
quotient singularities.
Therefore the rational cohomology groups of $\Y$ have \textit{pure}
Hodge structures. 
For the projective case, see \cite[Cor.\ 14.4]{Danilov}.
For the general case,
it follows from the following two facts: one is that
the intersection cohomology of every complete variety has pure Hodge 
structures \cite{Hodge-module}, the other is that the rational cohomology of
a variety with quotient  singularities equals the intersection
cohomology.


\subsection{The motivic measure on twisted $\infty$-jets}

Let
$\X$ be a smooth Deligne-Mumford stack of pure dimension $d$.
By abuse of notation, we also denote by $\overline { \cL_\infty \X }$
the set of points  $|\cL_\infty \X| = | \overline { \cL_\infty \X } |$.
We denote by $\pi_n$ the natural morphism 
$\overline {\cL_\infty\X} \to \overline {\cL_n \X}$.

\begin{defn}
A subset $A$ of $\overline {\cL_\infty \X} $
 is \textit{stable at level $n$} if we have:
\begin{enumerate}
	\item $\pi_n(A)$ is a constructible subset in $\overline {\cL_n \X}$,
        \item $A=\pi_n^{-1}\pi_n(A)$,
\end{enumerate}
A subset $A\subset \overline {\cL_\infty \X} $ is \textit{stable} if it is 
stable at level $n$ for some $n\in \Zz$.
\end{defn}

We define the notion of the \textit{measurable subset} similarly.
Then we define  
\textit{motivic measure} $\mu_\X$
 on  $\overline {\cL_\infty \X}$ by 
\[
  \mu_\X (A) := \LL^{-nd} \chi_h (\pi_n (A)) \in \hat K_0 (\HS), \ n \gg 0 , 
\]
where $\LL = \{ \QQ(-1)\} = \chi_h (\A^1) $. It is well-defined by the 
following:

\begin{prop}\label{piecewise}
Let $n\in \Zz$, let $B \subset \overline {\cL_n \X}$ be a constructible
subset and let $C$ be the inverse image of $B$ by the natural morphism
$\overline {\cL_{n+1} \X} \to \overline {\cL_{n} \X}$.
Then we have the equality
$\chi_h(C) = \LL^{d}\chi_h (B)$.
\end{prop}

\begin{proof}
Let $\Y$ be a connected component of $I^l(\X)(=\cL_0^l \X)$
and put $(\cL_n \X)_\Y:= \cL_n \X \times _{I(\X)} \Y$.
Let $ \gamma = (\tilde \gamma,(\overline \gamma, b), \id_{\overline \gamma})
\in (\cL_n \X)_\Y \subset L_m \X \times_\X \Y $ be an $\Omega$-point. 
Then we have the following commutative diagram of solid arrows:
\[\xymatrix{
   D_0 \otimes \Omega  \ar[r]^{(\overline \gamma, b)} 
\ar@{^{(}->}[d] & \Y \ar[r] &\cN \ar[d] \\
   D_m \otimes \Omega \ar[rr]_{\tilde \gamma} \ar@{-->}[urr]^{\sigma} &  &\X  .
}\]
Since $\cN \to \X$ is $\Y$-\'etale,
 there is a unique broken arrow $\sigma$ fitting into the diagram. 
Sending  $\gamma$ to  $\sigma$ determines a closed immersion
\[
\iota: (\cL_n \X)_\Y \hookrightarrow L_m \cN .
\]
Let $\mathfrak g$ be the canonical automorphism of $\cN$.
In view of the definition of
$\cN$ and $\clubsuit$ in the proof of Lemma \ref{closed}, we see 
that for $\sigma \in L_m \cN$, $\sigma \in \Img (\iota)$ iff
 $ \sigma \circ \zeta_l = \mathfrak g \circ \sigma $.
From Lemma \ref{linearization}, there exists
an atlas $h: V \to \Y$ such that
$\cN_V := \cN \times_\Y V \cong V \otimes  \CC[[v_1, \dots , v_c]]$, 
$(c = \dim \X - \dim \Y)$, and the pull-back of $\mathfrak g$ is given by
$\diag (\zeta_l^{a_1},\dots,\zeta_l^{a_c})$, $1 \le a_i \le l$.
 So for an $m$-jet $\delta$
 on $\cN_V$, $h \circ \delta $ is in $\Img (\iota)$ iff 
the image of $\delta$ by $L_m (\cN_V) \to L_m V$ is $\zeta_l$-invariant and
$\delta^*(v_i)$ is of the following form:
\[
   r_0 t^{a_i} + r_1 t^{a_i +l} + r_2 t^{a_i +2l} + \cdots .
\]
Therefore we have that  
$ (\cL_n \X)_\Y \times _\Y V \cong  L_n V \times \A^{nc}$ and
the projection 
$(\cL_{n+1} \X)_\Y \times _\Y V \to (\cL_n \X)_\Y \times _\Y V$ is a
Zariski locally
trivial $\A^d$-bundle. 

We may assume that a finite group, say $G$, acts on 
each connected component $V'$ of 
$V$, $V' \to \overline \Y$ is $G$-invariant and
the induced morphism $V'/G \to \overline \Y$ is \'{e}tale 
(see Lemma \ref{locallyquotient}).
Then we have that $\overline{(\cL_n \X)_\Y} \times _{\overline \Y} (V'/G)
\cong ((L_n V') /G) \times \A^{nc}$.
Therefore, from Lemma \ref{lem-hodgetrivial}, to prove the proposition,
it suffices to show that there is a stratification of $(L_n V') /G$
such that the natural morphism
$(L_{n+1} V') /G \to (L_n V') /G$ is, over each stratum,
 an analytically locally trivial fibration
of the quotient of an affine space by a linear
finite group action.

Let $H$ be a subgroup of $G$ and $W \subset V'$ a connected component of
the locus of the points with stabilizer $H$. 
Let $w \in W$ be a close point.
As is well known, 
there is a representation $\rho:H \subset \GL_{\dim \Y} (\CC)$ 
which describes 
the $H$-action on an analytic neighborhood of $w$.
Let $(L_n V)_{w} \subset L_n V$ be the subset of the jets which
maps the only point of $D_n$ to $w$.
Then the induced $H$-action on $(L_n V)_{w}
 \cong \A^{m \dim \Y}$ is given by 
$\rho ^{\oplus m} $.
Therefore $ (L_{n+1} V')_{W} /H \to (L_n V')_{W} /H $ is an
analytically locally trivial fibration of $\A ^{\dim \Y}/H$.
Let $G'\subset G$ be the subgroup of the elements
keeping $W$ stable. 
Then $H$ is a normal subgroup of $G'$.
It is easy to see that
the image of $(L_n V')_{W} /H$ in $(L_n V') /G$ is naturally isomorphic to
$((L_n V')_{W} /H)/(G'/H)$. Since $G'/H$ freely acts on $(L_n V')_{W} /H $,
the assertion follows.
\end{proof}

\begin{lem}\label{linearization}
Let $G$ be a finite group, $V$ a smooth $G$-variety and
$W$ a smooth closed subvariety consisting of
 $G$-invariant points.
\begin{enumerate}
	\item 
Assume that $V$ and $W$ are affine, say $V= \Spec R $ and $W= \Spec R/\fp$.
Moreover assume that 
$\fp$ is generated by $c = \codim (W,V)$ elements.
Then
the completion of $V$ along $W$ is
isomorphic as $G$-schemes to
$\Spec (R/\fp)[[x_1 , \dots, x_c]] $,  $G \subset \GL_c (R/\fp) $. 

  \item Assume that $G$ is a finite cyclic group.
Then there is an affine open covering $\cup V_i$ of $V$ such that
for every $i$, if we write 
$V_i = \Spec R$ and $W \cap V_i = \Spec R/\fp$, 
the completion of $V_i$  along $W \cap V_i$ is isomorphic as $G$-schemes to
$\Spec (R/\fp)[[x_1 , \dots, x_c]] $,  $G \subset \GL_c (\CC) $.
\end{enumerate}
\end{lem}

\begin{proof}
(1).
We denote by $\hat R$ the completion of $R$ with respect to an ideal
$\fp$ and by $\widehat{R_{\fp}}$ the completion of the local ring
$R_{\fp}$ with respect to the maximal ideal.
Let $K$ be the quotient field of $R/ \fp$ and 
$f_i$ generators of $\fp$.
It is well known that there is an isomorphism
$\widehat{R_{\fp}} \to K[[x_1, \dots, x_c]]$ sending $f_i$ to $x_i$.
We have a natural injection $\hat R \to \widehat{R_{\fp}}\cong
 K[[x_1, \dots, x_c]]$. Clearly the image contains the subring
$ (R/\fp)[[x_1, \dots, x_c]] $. 
Consider the injection
$\iota:(R/\fp)[[x_1, \dots, x_c]] \to \hat R$.
Since the induced map $R/ \fp \to \hat R / \hat \fp$ is the identity and
the images $f_i$ of $x_i$ generate $\hat \fp$, $\iota$ is a surjection
and hence an isomorphism (see \cite[Th. 7.16.]{Eisen}).

Consider the induced $G$-action on $\hat R = (R / \fp) [[x_1 , \dots, x_c]]$.
For $g\in G$, write $g(x_i) = \sum a_{ij} x_j + (\text{higher terms})$,
$a_{ij}\in R / \fp$. Denote by $\bar g$ the endomorphism of $\hat R$
associated to the invertible matrix $(a_{ij})$.
Working in characteristic zero, we can isomorphically replace $x_i$ with 
$x'_i = \sum _{g \in G} \bar g ^{-1} g (x_i) $. We find that
with respect to the new
coordinates, the $G$-action is linear.

(2). In the last situation,
$\Spec \hat R$ is naturally isomorphic to the completion of 
the normal bundle $N_{W/V}$ along the zero section.
Here let us assume $G$ is a cyclic group with generator $g$. 
Then $N_{W/V}$ decomposes to eigenbundles.
On each eigenbundle, the $g$-action is uniquely represented by a scalar
matrix $aI$ where $I$ is the identity matrix and $a \in \CC$.
Therefore shrinking $V$ to an open subset where the eigenbundles are free,
we conclude that the $G$-actions on $N_{W/V}$ and $\Spec \hat R$ are realizable
in $\CC$.
\end{proof}

\begin{rem}
The author guesses that even in the case of a general finite group,
the action on
$N_{W/V}$ is \textit{\'{e}tale} locally realizable in $\CC$.
From facts on splitting fields of finite groups 
(see \cite{curtis-reiner}), this is true at least over the generic point of 
$W$.
\end{rem}

\begin{lem}\label{lem-hodgetrivial}
Let $T$ and $S$ be varieties and $f : T \to S$ an analytically locally trivial
fibration of $\A ^d/G$ for a finite group $G \subset \GL_d (\CC)$.
Then $\chi _h (T) = \chi _h (S) \LL ^d$.
\end{lem}

\begin{proof}
Since the fiber is a quotient of an affine space, the higher direct images of
$\QQ_T$ vanishes;
\[
R^i f_* \QQ _T \cong \begin{cases}
  \QQ _S & ( i = 0) \\
  0 & (i > 0) .
\end{cases}
\]
Hence the spectral sequence is degenerate and it follows that 
$H^i ( T, \QQ ) \cong H ^i (S , \QQ)$ for every $i$.

Taking a stratification of $S$, we may assume that $S$ is smooth.
Since $S$ and $T$ have at most quotient singularities (in the analytic sense), 
by Poincar\'{e} duality, we conclude that $H ^{\dim T - i} _c( T , \QQ)
\cong H^ {\dim S - i} _c(S , \QQ)$. 

Regarding the sheaves as mixed Hodge modules, studied by Saito 
\cite{Hodge-module} (see also \cite{Intro-Hodge}), we can
regard the isomorphisms above of cohomology
groups as ones of mixed Hodge structures.
This implies the assertion.
\end{proof}


\subsection{Main theorem}

Let $X$ be a variety with Gorenstein quotient singularities.
Then $X$ has canonical singularities.
Let $\X$ be a smooth Deligne-Mumford stack without reflections such that
$X$ is the coarse moduli space of $\X$.
We denote by $\lambda$ the canonical morphism $\X \to X$.
If $\gamma:\D_n^l\otimes \Omega \to \X$ is
 a twisted $n$-jet on $\X$ of order $l$,
then it induces a morphism $\gamma':D_n\otimes \Omega \to X$
of the coarse moduli spaces.
We define the map $\lambda_{(n)} : \overline{\cL_n \X} \to L_n X$ by 
$\gamma \mapsto \gamma'$.
The following is our main result.

\begin{thm}\label{main}
Let $B \subset \overline {\cL_\infty \X}$ be a measurable subset 
and put $A:= \lambda_{(\infty)}(B)$.
Then we have the following equation in $\hat K_0 (\HS)$:
\[
\chi _h \muGor_X (A) = 
\sum_{\Y \subset I(\X)} \LL^{s(\Y)}\mu_\X 
(\pi_0^{-1}(\overline \Y) \cap B ),
\]
where $\Y$ runs over the connected components of $I(\X)$.
\end{thm}

The proof is postponed until the end of the section.

\begin{cor}\label{cor-of-main}
Let $X$ and $X'$ be complete  varieties with Gorenstein quotient 
singularities. Suppose that there are proper birational morphisms $Z\to X$ and
$Z\to X'$ such that $K_{Z/X}=K_{Z/X'}$.
Then the orbifold cohomology groups
 of $X$ and $X'$ have the same Hodge structure.
\end{cor}

\begin{proof}
By Theorem \ref{main} and Proposition \ref{piecewise}, we have
\begin{align*}
\chi_h \muGor_X (L_\infty X) &= 
\sum_{\Y \subset I(\X)} \LL^{s(\Y)} \chi_h(\overline \Y) \\
& = \sum_i (-1)^i \{\Horb^i(X,\QQ)\}.
\end{align*}
From Lemma \ref{Goren}, we have 
$\sum_i (-1)^i \{\Horb^i(X,\QQ)\} = \sum_i (-1)^i \{\Horb^i(X',\QQ)\}$.
Since $\Horb^i(X,\QQ)$ and $\Horb^i(X',\QQ)$ have a pure Hodge 
structure of weight $i$, 
 $\{\Horb^i(X,\QQ)\}=\{\Horb^i(X',\QQ)\}$ for every $i$.
\end{proof}

\begin{lem}\label{bijective}
Let $X_{\sing}$ denote the singular locus of $X$ with reduced subscheme
structure.
\begin{enumerate}
	\item The subset $\lambda _{(\infty)}^{-1}(L_\infty ( X_\sing ))$ 
             of $ \overline {\cL_ \infty \X} $ is of
             measure zero.
        \item The map $\lambda _{(\infty)}$ is bijective over 
          $L_\infty X \setminus L_\infty (X_{\sing})$. 
\end{enumerate}
\end{lem}

\begin{proof}
(1). It suffices to show that for every $n$, 
$\pi_n(\overline {\cL_\infty \X} \setminus
\lambda _{(\infty)}^{-1}(L_\infty ( X_\sing ))) = \overline{\cL_n \X} $.
But this is clear by the local description of $\overline{\cL_n \X}$ in 
the proof of Proposition \ref{piecewise}.

(2). Surjectivity. Let $\eta:\Spec \Omega[[t]] \to X$ be an $\Omega$-point of 
$L_\infty X \setminus L_\infty ( X_{\sing} )$ with an algebraically
closed field $\Omega$.
We define $\D$ to be the normalization of 
the fiber product $\X \times_X \Spec \Omega[[t]]$, (see Definition 
\ref{normalization}).
 Then the Deligne-Mumford stack $\D$ contains
the scheme $\Spec \Omega ((t))$ as open substack. Therefore
the coarse moduli space 
$ \overline \D $ of $\D$ contains $\Spec \Omega ((t))$ as
open subscheme. The scheme $ \overline \D $ 
must be the spectrum of a local ring
$S \subset \Omega((t))$. From the universality of coarse moduli space, there
is a natural morphism $  \overline \D \to \Spec \Omega[[t]] $. So we have 
$\Omega[[t]] \subset S \subset \Omega((t))$. Because 
$\Omega[[t]]$ and $\Omega((t))$ are 
the only intermediate rings between $\Omega[[t]]$ and $\Omega((t))$,
the ring $S$ must be $\Omega[[t]]$.
Suppose that $E = \Spec \tilde S$ is an atlas of $\D$ and 
$\tilde S$ is a regular local ring. Since $S=\Omega[[t]]$ is henselian, 
the natural morphism $ E \to  \overline \D  $ is finite 
(\cite[IV. Th.\ 18.5.11]{EGA}). 
Hence $\tilde S$ is complete, (see \cite[Cor.\ 7.6]{Eisen}). 
So $\tilde S \cong \Omega [[t]]$. 
Consider the groupoid space 
$E \times_\D E \rightrightarrows E$.
The scheme  $E \times_\D E$
 must be the disjoint sum of spectra of
complete regular local rings. 
Since the first projection
$\pr_1:E \times_\D E \to E$ is \'etale,
there is an isomorphism 
\[
   E \coprod \dots \coprod E \cong 
   E \times_\D E
\]
such that the composition 
\[
E \coprod \dots \coprod E \cong
   E \times_\D E \xrightarrow { \pr_1 } E
\]
is isomorphic on each component.
If  $l$ denote the number of the components in $ E \times_\D E $,
then the second projection $\pr_2:E \times_\D E \to E$ 
determines the action of some group $G$ on $E$ with $ |G| = l $.
Since this action is effective, the group $G$ is isomorphic to $\mu_l$
for some $ l $.
 For a suitable isomorphism $\mu_l \cong G$,
the action is given by $t \mapsto \zeta_l t$.
Hence the stack  $\D$
is isomorphic to $\D_\infty^l \otimes \Omega$ and the morphism 
$\D_\infty^l \otimes \Omega \cong \D \to \X$ is a twisted $\infty$-jet on $\X$.
The image of this twisted $\infty$-jet by $ \lambda _{(\infty)}$ is $\eta$.

Injectivity. Let $\gamma _1 , \gamma _2 :\D_\infty^l\otimes \Omega \to \X$ be
 two twisted $\infty$-jets on $\X$ of order $l$. 
We suppose that $\eta := \lambda _{(\infty)} (\gamma _1 )
  = \lambda _{(\infty)} (\gamma _1 )$ and 
  $\eta \in L_\infty X \setminus L_\infty ( X_{\sing} )$. 
Construct $\D$ from $\eta$ as above. Then 
for each $i \in \{ 1,2 \}$, there is a unique morphism 
$ h_i : \D_\infty^l\otimes \Omega \to \D$ such that the following diagram
is commutative:
\[
 \xymatrix{
 \Spec \Omega ((t)) \ar@{=}[d] \ar@{^{(}->}[r] & 
    \D_\infty^l\otimes \Omega \ar[d] ^{h_i} \ar[dr] ^{\gamma _i} &  \\
 \Spec \Omega ((t)) \ar@{=}[d] \ar@{^{(}->}[r] & \D \ar[d] \ar[r]
          & \X \ar[d] ^{\lambda}   \\
 \Spec \Omega ((t)) \ar@{^{(}->}[r] & D_\infty \otimes \Omega
 \ar[r] ^{\eta} & X .
 }
\]
Let $E$ be an atlas of $\D$ as above. Then the natural morphism
$E \times_{ \D, h_i } (\D_\infty^l\otimes \Omega) \to E$ is a birational 
morphism of smooth 1-dimensional schemes. Therefore it is an isomorphism and 
so is $h_i$ (\cite[Prop.\ 3.8.1]{LMB}).
Then we can easily see that $\gamma_1$ and $\gamma_2$ have the same image
in $ | L_\infty \X \times _\X I^l(\X) |  $.
\end{proof}

To prove Theorem \ref{main},
 we now need to generalize the transformation rule.
Let $\cV$ be a Deligne-Mumford stack over $D_\infty$ of 
pure relative dimension $d$.
For each $n \in \Zz \cup \{ \infty \}$, we define 
$\cV_n$ to be the moduli stack of the $D_\infty$-morphisms $D_n \to \cV$.
Then for $m \ge n$, there is a natural projection $\cV_m \to \cV_n$.
So we can define the \textit{motivic measure}
 $\mu_\cV$ over $\overline{\cV_\infty}$ 
which takes values in $\hat K_0 (\HS)$,
in a similar fashion as before. (We should replace condition (3) in 
Definition \ref{def-stable} with the condition that 
$\chi _h (\pi_{m+1} A) = \LL ^d \chi _h (\pi_m A)$.
It makes sense because of Lemma \ref{lem-hodgetrivial}.)

Let $\cW$ be another Deligne-Mumford stack over $D_\infty$ of
pure relative dimension $d$ and let $h : \cW \to \cV$ be a 
$D_\infty$-morphism.
We put $\Omega'_{\cW/\cV}:= 
\Img ( \Omega_{\cW/ D_\infty } \setminus (tors)  \to \Omega _{\cW / \cV} ) $,
where $ (tors) \subset \Omega_{\cW/ D_\infty } $ 
is the  torsion.
Then we define the \textit{jacobian
 ideal sheaf} $\J _ h$ of $h$ to be the $0$-th Fitting
ideal of $ \Omega'_{\cW/\cV}$.

\begin{thm}\label{thm-GeneralTrans}
Let $A \subset \overline{\cW_\infty}$ be a measurable set.
Suppose that $h_\infty|_A : A \to \overline{\cV_\infty}$ is injective.
 Let $\nu$ be a measurable function
on $h_\infty(A)$. Then
\[
\int_{h_\infty (A)} \LL^\nu  d\mu_\cV = 
  \int_A \LL^{\nu \circ h_\infty - \ord \, \J_h} d\mu_{\cW}.
\]
\end{thm}

\begin{proof}
It is a direct consequence of the following lemma.
\end{proof}

We denote by $\J(\cV/D_\infty)$ (resp. $\J(\cW/D_\infty)$) the
$d$-th Fitting ideal sheaf of $\Omega_{\cV / D_\infty}$ 
(resp. $\Omega_{\cW / D_\infty}$).

\begin{lem}\label{lem-key}
Let $ A \subset \overline{\cW_\infty} $ be a stable subset of 
level $l$. Assume that $h_{\infty} |_A$ is injective, that 
$\ord \J _ h$ is constant equal to $e < \infty$ and
 that $ \ord \J(\cV/D_\infty)$  and  $\ord \J(\cW/D_\infty)$ are
 bounded from above on $h_\infty (A)$ and  $A$ respectively.
Then for
 $n \gg 0$, $ h_n : \pi _n A \to h_n \pi_n A$ is a piecewise
trivial $\A^e$-bundle.
\end{lem}

\begin{proof}
Looijenga's proof \cite[Lem.\ 9.2]{looi} works also in this setting.

Take a non-twisted $\infty$-jet $\gamma:\Spec \Omega[[t]] \to \cW $ in $A$,
 and put $\fm := (t) \subset \Omega[[t]]$. 
Let $q$ be the image of the closed point by $\gamma$.
Take another $\theta \in A$ such that 
$ \pi_{n-e}  (\gamma) =  \pi_{n-e}  (\theta) $.
Then the morphism 
\[
 \theta^* - \gamma^* : \cO_{\cW, q} \to 
 \fm ^{n-e+1} / \fm ^{2(n-e+1)} 
\]
 is a $\CC[[t]]$-derivation. 
So it defines an $\Omega[[t]]$-module homomorphism 
\[
 \partial \theta   :  \gamma^* \Omega_{\cW / D_\infty} 
\to \fm^{n-e+1}/\fm^{2(n-e+1)}. 
\] 
The length of the torsion of $\gamma^* \Omega_{\cW / D_\infty}$ 
equals $(\ord \J(\cW / D_\infty)) (\gamma)$, 
hence it is bounded. So, since $n \gg 0$,
 the composition map 
\[
 \overline{ \partial \theta } : \gamma^* \Omega_{\cW / D_\infty} 
\xrightarrow{ \partial \theta } \fm^{n-e+1}/\fm^{2(n-e+1)} \to
\fm^{n-e+1}/\fm^{n+1}
\]
 annihilates the torsion.
Conversely, every $\Omega[[t]]$-module homomorphism
$\gamma^* \Omega_{\cW / D_\infty} \to \fm^{n-e+1}/\fm^{n+1}$
which annihilates the torsion is $\overline{ \partial \theta } $
for some $\theta$.

After some works, we can see that
if $\theta \in A$ is such that $h_n \pi_n  (\gamma) = 
h_n \pi_n  (\theta) $,
then $ \pi_{n-e}  (\gamma) =  \pi_{n-e}  (\theta) $, (see
\cite[Lem.\ 9.2]{looi}). So $\overline {\partial \theta}$ is defined.
It is easy to see that $\pi_n(\theta_1) = \pi_n(\theta_2)$ iff
$ \overline{\partial \theta_1} =\overline{ \partial \theta_2} $ , and that 
$h_n \pi_n(\theta_1) =h_n \pi_n(\theta_2) $ iff 
$ \overline{\partial \theta_1} $ and $\overline{ \partial \theta_2} $
have the same image in 
$\Hom_{\Omega[[t]]}( (h\gamma)^* \Omega_{\cV / D_\infty},
\fm^{n-e+1}/\fm^{n+1})$.
Hence  $ h_n^{-1} h_n \pi_n (\gamma)  $ is isomorphic to an
affine space, 
\begin{align*}
\Hom_{\Omega[[t]]} ( \gamma^* \Omega'_{\cW / \cV} & ,
\fm^{n-e+1}/\fm^{n+1}) \\
\cong
\Ker & (
\Hom_{\Omega[[t]]} ( \gamma^* \Omega_{\cW / D_\infty}\setminus (tors),
\fm^{n-e+1}/\fm^{n+1}) \\
 & \to
\Hom_{\Omega[[t]]}( (h\gamma)^* \Omega_{\cV / D_\infty},\fm^{n-e+1}/\fm^{n+1})
).
\end{align*}
The length of $\gamma^* \Omega'_{\cW / \cV}$ equals $e = \ord \J _ h (\gamma)$.
So $ h_n^{-1} h_n \pi_n (\gamma)$ is isomorphic to an affine space of 
dimension $e$. 

The rest is easy.
\end{proof}

\begin{proof}[Proof of Theorem \ref{main} ]
Let $\Y$ be a connected component of $I^l(\X)$ and 
$\cN$ its formal neighborhood. 
We may assume that $B$ is 
 contained in $\pi_0^{-1}(\overline \Y)$.
Let $\tilde \cN$ be the quotient of $\cN$ by 
the canonical automorphism $\mathfrak g$, that is, 
$\cSpec \cA^{\mathfrak g}$ where $\cA^{\mathfrak g} \subset \cA$ is the 
subsheaf of the $\mathfrak g$-invariant sections.
Then the natural morphism $\cN \to X$ factors as 
\[
\cN \to \tilde \cN \xrightarrow{f} X.
\]

In the proof of Proposition \ref{piecewise}, we saw that 
for each $n,\, m$ with $m=nl$,
there is a closed immersion
$\iota:(\cL_n \X)_\Y \hookrightarrow L_m \cN$.
Let $[l]:D_m \to D_n$ be the morphism associated to
the ring homomorphism defined by $t \mapsto t^l$.
If $\gamma: \D_n^l \otimes \Omega \to \X$ is a twisted $n$-jet in
$\pi_n (B)$, then $\iota(\gamma)$ fits into the diagram 
\begin{equation}\label{diag1}
\xymatrix{
D_m \otimes \Omega \ar[r]^{\iota(\gamma)} \ar[d]_{[l]} & \cN \ar[d] \\
D_n \otimes \Omega \ar[r]_{\sigma} &  \tilde \cN .
}
\end{equation}
Then $\lambda_{(n)} (\gamma) = f_n ( \sigma) $. 
We define a subset 
$\tilde B \subset \overline {L_\infty \tilde \cN}$ 
to be the image of $B$ by the map $\gamma \mapsto \sigma$.
Then $A = f_\infty (\tilde B) $ and $f_\infty|_{\tilde B}$ is bijective
outside of subsets of measure zero. 
Let $\I_X$ (resp. $\I_{\tilde \cN} $) be the ideal sheaf on $X$ defined by
\begin{align*}
   \I_X \omega_X &= \Img ( \Omega_X^d  \to \omega_X) \\
  (\text{resp.}\ \I_{\tilde \cN} \omega_{\tilde \cN} &= 
  \Img ( \Omega_{\tilde \cN}^d  \to \omega_{\tilde \cN}) ),
\end{align*}
and define $\muGor_X$ and $\muGor_{\tilde \cN}$ 
to be $\LL^{\ord \I_X} \mu_X$ and 
$\LL^{\ord \I_{\tilde \cN}} \mu_{\tilde \cN}$ respectively.
Since the morphism $f$ has no ramification divisor, 
by a similar argument as the proof of Lemma \ref{Goren},
we see
$ f^{-1} \I_X  =  \J_f \cdot \I_{\tilde \cN} $, 
where $\J_f$ is the jacobian ideal sheaf.
So, by Theorem \ref{thm-GeneralTrans}, we obtain 
\begin{equation*}\label{tototo}
 \muGor_{\tilde \cN}(\tilde B) = \muGor_{X}(A) .
\end{equation*}
We have thus reduced the problem to the case of a cyclic quotient;
it suffices to show the following lemma.
\end{proof}

\begin{lem}
Let the notation as above.
We have $ \LL ^{s(\Y)} \mu _ \X (B)
= \chi_h \muGor _{\tilde \cN} (\tilde B)$.
\end{lem}

\begin{proof}
The proof is essentially by a trick used in \cite{DL-quotient}.
We first consider an easy case where $\X$ is a quotient stack 
$[ \A^c_R / G]$ of
an affine space over a ring $R$ whose spectrum is a smooth variety of
dimension $d-c$,
 and $G \subset \SL _c (\CC)$ a finite cyclic group of order $l$ generated by 
$g= \diag (\zeta_l^{a_1},\dots,\zeta_l^{a_c})$, $1 \le a_i < l$.
Suppose that $\Y$ is the component associated to $g$.
Then $\cN = \hat \A ^c_R (=\Spec R [[x_1 , \dots , x_c]])$, its canonical
automorphism is $g= 
\diag (\zeta_l^{a_1},\dots,\zeta_l^{a_c})$, 
and $\tilde \cN = \hat \A ^c _R/ G$.
Since the natural morphisms $\hat \A ^c_R \to \A ^c_R$ and
$\hat \A ^c_R/G \to \A ^c_R/G$ are $(\Spec R)$-\'{e}tale, and since
we consider only  jets which send the only closed point into $\Spec R$,
 it makes no matter
to replace $\hat \A ^c_R$, $\hat \A ^c_R / G$ with
$ \A ^c_R$, $ \A ^c _R/ G$.

Consider three $R$-algebra 
homomorphisms; (i)  $u^* : R[[t]][\underline{x}] \to 
R[[t]][\underline{x}]$, $x_i \mapsto t^{a_i} x_i$, 
(ii) $\alpha^* : R[[t]][\underline{x}] \to R[[t]][\underline{x}]$, 
$ x_i \mapsto x_i$, $t \mapsto t^l $,  (iii) $\beta^* : 
R[[t]][\underline{x}]^{G} \to R[[t]][\underline{x}]$, the composition
of $\alpha^*$ and the inclusion $R[[t]][\underline{x}]^{G} \hookrightarrow 
R[[t]][\underline{x}]$.
Since $R[\underline {x}]^{G}$ is generated by
the monomials $x_1^{m_1} \dots x_c ^{m_c} $ with $\sum a_i m_i \equiv 0 \
(\mathrm{mod} \ l)$, there is a $R[[t]]$-homomorphism 
$v^*: R[[t]][\underline{x}]^{G} \to R[[t]][\underline{x}] $ with
$u^* \circ \beta^* = \alpha^* \circ v^*$;
\[\xymatrix{
R[[t]][\underline{x}] \ar[r]^{u^*}  & R[[t]][\underline{x}] \\
R[[t]][\underline{x}] ^{G}\ar[r] _{v^*} \ar[u]^{\beta^*} & R[[t]][\underline{x}] \ar[u] _{\alpha^*} .
}
\]
Here the horizontal arrows are $R[[t]]$-algebra 
homomorphisms and the vertical ones
send $t \mapsto t^l$. 
Write the diagram of the associated schemes as follows;
\[\xymatrix{
\cN [[t]] \ar[d]_{\beta}  & E_2 [[t]] \ar[d] ^{\alpha} \ar[l] _{u}\\
\tilde \cN [[t]] & E_1 [[t]] \ar[l]^{v} .
}
\]
Here $E_i$ are copies of $\A^c_R$.

Let $B_0$ be the image of $B$ by $\iota:(\cL_{\infty} \X)_\Y 
\hookrightarrow L_{\infty} \cN$.
Then for $\gamma' \in B_0$,
$(\gamma ')^* (x_i) $ is of the form,
\begin{equation}\label{eq-temp}
(\gamma ')^* (x_i) = 
 r_0 t^{a_i} + r_1 t^{a_i +l} + r_2 t^{a_i +2l} +
\cdots  .
\end{equation}
If we put $ \eta := u_{\infty}^{-1} (\gamma') $, then we have
\[
 \eta ^* (x_i) = 
 r_0  + r_1 t^{l} + r_2 t^{2l} +
\cdots  .
\]
Therefore if we define $\delta \in L_{\infty} E_1 $ by 
\begin{equation}\label{eq-temp2}
\delta ^*(x_i) = 
 r_0  + r_1 t + r_2 t^{2} + \cdots,
\end{equation}
then  we have the following commutative diagram,
\[\xymatrix{
E_2[[t]]  \ar[d]_{\alpha} &
 D_{\infty} \otimes \Omega \ar[d]^{[l]} \ar[l]_{\eta}\\
E_1[[t]]  & D_{\infty} \otimes \Omega \ar[l]^{\delta}.
}
\]
 Here $[l]$ is the morphism
defined by $t \mapsto t^l$.
Let $B_1 \subset L_{\infty}(E_1)$ be the image of $B_0$ by the map
$\gamma' \mapsto \delta$.
It is easy to see that
if  $B$ is stable at level $n$, then so is
 $B_1$, and that
$ \{\pi_n(B)\} =\{ \pi _{nl}( B_0)\} = \{\pi_{n} (B_1)\} \LL^{-c}$,
where $\pi_i$ are truncation morphisms of $\overline{\cL_{\infty} \X}$, 
$L_{\infty} \cN$, and $L_{\infty } E_1$ respectively.
Therefore 
\begin{equation}\label{BB2}
 \chi _ h \mu _ {E_1} (B_1) =
 \mu_{\X}(B) \LL^{c} .
\end{equation}

Put $\sigma := v_{\infty} (\delta)$. The chain of the correspondences,
$\gamma \mapsto \gamma'  \mapsto
\delta \mapsto \sigma$, defines a map $(\cL_{\infty} \X)_\Y \to 
L_{\infty} \tilde N $,
which is the same as one in the proof of the theorem 
(see diagram (\ref{diag1}) and compare it with the last two ones).

Shrinking $\Spec R$ to an open subset, suppose that the canonical sheaf 
$\omega _{\Spec R}$ of  $\Spec R$ is generated by
a section $e'$.
Consider a $d$-form $e= dx_1 \wedge \dots \wedge dx_c \wedge e'$  on $\cN$.
This is stable under the $G$-action.
If $r$ denotes the natural morphism $\tilde \cN \to \cN$,
the canonical sheaf $\omega _{\tilde \cN}$ of $\tilde \cN$
is generated by a $d$-form $\tilde e$ with $ r^{*} \tilde e=e$. 
Direct computation 
gives $ v ^ * \tilde e = t^{\sum a_i /l} 
(dx_1 \wedge \dots \wedge dx_c \wedge e')$.
Hence we have the following equations of subsheaves of 
$\Omega^d_{E_1[[t]]/D_{\infty}}$,
\begin{align*}
(t^{\sum a_i /l})v^{-1}\J _{\tilde \cN} \cdot 
v^* \omega_{\tilde \cN/D_{\infty}}
&= (t^{\sum a_i /l}) (v^* \Omega^d_{\tilde \cN[[t]]/D_{\infty}})/(tors) \\
&= (t^{\sum a_i /l})\J _{v} \cdot \Omega^d _{E_1[[t]]/D_{\infty}} \\
&= \J _{v} \cdot v^* \omega_{\tilde \cN/D_{\infty}}.
\end{align*}
This means that $\ord \J _{\tilde \cN} \circ v_{\infty} -
\ord \J_{v} \equiv - \sum a_i /l$.
From the transformation rule, we obtain that
$\chi_h \muGor _ {\tilde \cN} (A) = \LL ^ {-\sum a_i /l} \mu _{E_1} (B_1) $, 
and
granting (\ref{BB2}), that  $\chi_h 
\muGor _ {\tilde \cN} (A) = \LL ^{s(\Y)} \mu_{\X} (B) $.
We have proved the assertion in this case.

As for the general case, the proof follows along almost the same lines:
We take the fiber product $ \cN \times _{ \Y} V$ for 
an atlas $V \to \Y$ to linearize the canonical automorphism. 
Define $B_V$, $B_{V,1}$ and $\tilde B_V$ in the evident fashion.
By the  argument for the preceding case and a common argument, 
after replacing $B$, 
we obtain that $B_{V,1}$ and $\tilde B_V$
are stable at level $n$ and  a morphism $v_n:
\pi_n B_{V,1} \to \pi_n \tilde B_V$ 
is a 
 trivial affine space bundle of the expected relative dimension.  
Here we have used Lemma \ref{lem-key} instead of
the transformation rule itself.
The natural morphism
$\pi _{n+1} B_V \to \pi_n B_{V,1}$ is an affine space bundle of
relative dimension $d-c$ which is trivial Zariski locally on $V$.
(Recall (\ref{eq-temp}) and (\ref{eq-temp2}). This bundle results from 
the truncation $L_{n+1} (\Spec R) \to L_{n} (\Spec R) $ and the identity of
$(r_i)_{0 \le i \le n} $.)
Hence $\pi _{n+1} B_V \to \pi_n \tilde B_V$ is also a Zariski locally
trivial affine space bundle of the expected relative dimension.
By the same argument as the proof of Lemma \ref{piecewise}, 
we can conclude that $\pi _{n+1} B \to \pi_n \tilde B$ is an
analytically locally trivial fibration of a quotient of
an affine space.
The lemma follows from
 Lemma \ref{lem-hodgetrivial}.
\end{proof}


\section{General results on Deligne-Mumford stacks}

In this section, we give some general results on Deligne-Mumford stacks
which we need in the preceding section.
There are good references for stacks today
(for example \cite{DM}, \cite{Vistoli}, \cite{Gomez} and \cite{LMB}).

We fix a base scheme $S$.

\subsection{Deligne-Mumford stacks}\label{subsec-DM}

A \textit{stack} is a category fibered in groupoids over $\schS$ 
such that every Isom functor is a sheaf and every descent datum is effective.
A morphism of stacks $\X\to\Y$ is \textit{representable}\footnote{In 
\cite{LMB}, this is called \textit{sch\'ematique}.
} if for any 
$U\in\schS$ and any morphism $U\to\Y$, the fiber product $U\times_\Y \X $
is represented by a scheme.

\begin{defn}
Let  \textbf{P} be a property of morphisms $f:Y\to X$ of $S$-schemes, 
stable under base change and local in the \'etale topology on $X$
(for example: surjective, proper \textit{etc}).
We say that
 a representable morphism $f:\Y \to \X$ of stacks has 
 property \textbf{P} if for every $S$-scheme $U$ and every morphism
$U \to \X $, the projection $U \times_\X \Y \to U$ has property \textbf{P}.
\end{defn}

\begin{defn}\label{DMstack}
A \textit{(separated) Deligne-Mumford stack} is a stack 
$\X$ which satisfies the
following:
\begin{enumerate}
	\item the diagonal $\Delta: \X \to \X\times\X$ is representable
          and finite,
        \item there exists a scheme $M$ and a 
          morphism $M\to\X$ (necessarily representable after (1)),
          which is \'etale and surjective.
\end{enumerate}
A scheme $M$ in (2) is called an \textit{atlas} of $\X$.
A Deligne-Mumford stack $\X$ is 
\textit{of finite type} 
if there is an atlas of finite type.
\end{defn}

\begin{defn}
Let \textbf{P} be a property of morphisms $f:Y\to X$ of $S$-schemes, 
stable under \textit{\'etale}
 base change and local in the \'etale topology on $X$
(for example: birational, being an 
open immersion with dense image \textit{etc}).
We say that
a representable morphism $f:\Y \to \X$ of stacks has 
property \textbf{P} if for every scheme $U$ and every \'etale morphism
$U \to \X $, the projection $U \times_\X \Y \to U$ has property 
\textbf{P}.
\end{defn}

\begin{defn}
Let \textbf{Q} be a property of schemes local in the \'etale topology
(for example: reduced, 
smooth, normal, locally integral
\textit{etc}).
Let $\X$ be a Deligne-Mumford stack.
We say that $\X$ has  property \textbf{Q} if an atlas of $\X$ has property
\textbf{Q}.
\end{defn}

\begin{defn}
A (not necessarily representable) morphism $f:\Y \to \X$ of 
Deligne-Mumford stacks of finite type is \textit{proper}
if there is a $S$-scheme $Z$ and
a proper surjective morphism $g:Z \to \Y$ such that $f\circ g$ is 
(necessarily representable and) proper.

A Deligne-Mumford stack $\X$ of finite type is \textit{complete} if it
is proper over $S$.
\end{defn}

Although our condition appears weaker than one of \cite[Def.\ 4.11]{DM}, 
they are actually equivalent by Chow's lemma 
(see \cite[Def.\ 4.12]{DM}, \cite[Th.\ 16.6]{LMB}, 
\cite[Prop.\ 2.6]{Vistoli}).

\begin{expl}
Let $Z$ be a $S$-scheme and $G$ a finite
group  acting on $X$. The \textit{quotient
stack} $[Z/G]$ is defined as follows: an object over $U\in\schS$ is a
$G$-torsor $P\to U$ with a $G$-equivariant morphism $P\to Z$, and  a
morphism over $U' \to U$ is a cartesian diagram
\[ \xymatrix{ \ar @{} [dr] |{\square}
P'\ar[r] \ar[d] & P \ar[d] \\
U' \ar[r] &  U
} \]
which is compatible with the $G$-equivariant morphisms $P'\to Z$ and $P\to Z$.
It is a Deligne-Mumford stack with a canonical atlas $Z \to [Z/G]$.
\end{expl}

Here we define points of a Deligne-Mumford stack. 
For details, see \cite[Ch.\ 5]{LMB}.

\begin{defn}
Let $\X$ be a Deligne-Mumford stack. A \textit{point} of $\X$ is 
a $S$-morphism $\Spec K \to \X$ for a field $K$ with a morphism
$\Spec K \to S$.
\end{defn}

Let $x_i : \Spec K_i \to \X \, (i=1,2)$ be points of $\X$. 
We say that $x_1$ and $x_2$ are \textit{equivalent} if 
there is a field $K_3$ such that $K_3 \supset K_1,K_2$ and 
the diagram
\[\xymatrix{
{\Spec K_3} \ar[r] \ar[d] & {\Spec K_2} \ar[d] \\
{\Spec K_1} \ar[r]        &   {\X} 
}\] 
commutes.

\begin{defn}
We define the \textit{set of points} of $\X$, denoted $|\X|$,
 to be the set of the equivalence classes of points of $\X$. 
\end{defn}

The \textit{Zariski topology} on $|\X|$ is defined 
as follows: an open subset is $|\U| \subset |\X|$ for an open substack
$\U \subset \X$. There is a 1-1 correspondence between the closed subsets
of $|\X|$ and the reduced closed substacks of $\X$.

We now introduce the notion of (\'etale) groupoid space which is 
equivalent to Deligne-Mumford stack. Some references are 
\cite[p.\ 668]{Vistoli}, \cite[(2.4.3), (3.4.3), Prop.\ 3.8, (4.3)]{LMB} and
\cite[Subsec.\ 2.4]{Gomez}.

\begin{defn}
An \textit{(\'etale) groupoid space} consists of the following data:
\begin{enumerate}
	\item two $S$-schemes $X_0$ and $X_1$,
      \item five morphisms: \textit{source} and \textit{target}
            $q_i: X_1 \to X_0 \; (i=1,2)$,
            \textit{origin} $\varepsilon: X_0 \to X_1$,
            \textit{inverse} $\tau : X_1 \to X_1$ and 
            \textit{composition} $m: X_1 \times_{q_1,X_0,q_2} X_1 \to X_1$ 
            which satisfies the following:
\begin{enumerate}
	\item $q_1$ and $q_2$ are \'etale and 
        $(q_1,q_2): X_1 \to X_0 \times X_0$ is finite,
        \item the axioms of \textit{associativity}, \textit{identity element} 
and
\textit{inverse}.
\end{enumerate}
\end{enumerate}
We denote this groupoid space by $X_1 \rightrightarrows X_0$.
\end{defn}

Given a groupoid space $X_1 \rightrightarrows X_0$, we define the category
fibered in groupoids
$[X_1 \rightrightarrows X_0]'$ as follows: an object over 
$U\in\schS$ is a morphism
$U \to X_0$ of $S$-schemes and a morphism of $a:U \to X_0$ to $b:V \to X_0$ is
a pair of morphisms $f:U \to V $ and $h : V \to X_1 $ such that 
$q_1 \circ h \circ f = a $ and $ q_2 \circ h = b$.
Then $[X_1 \rightrightarrows X_0]'$ is a \textit{prestack}
(see \cite[3.1]{LMB}).

\begin{defn}
We define the \textit{stack associated}  to a groupoid space
$X_1 \rightrightarrows X_0$, denoted $[X_1 \rightrightarrows X_0]$,
 to be the stack associated to the prestack
$[X_1 \rightrightarrows X_0]'$ (\cite[Lem.\ 3.2]{LMB}). 
\end{defn}

The stack $\X=[X_1 \rightrightarrows X_0]$ is a Deligne-Mumford stack with
a canonical atlas $X_0 \to \X$.
We can identify the fiber product $X_0 \times_{\X} X_0$ with $X_1$.
Conversely, given a Deligne-Mumford stack $\X$ and an atlas $X_0 \to \X$,
then the schemes $X_0$ and $X_1 =M \times_\X M$
 underlies a natural groupoid space structure
with $q_i= \pr_i$ and $\varepsilon = \Delta_{M/\X}$.
The associated stack $[M \times_\X M \rightrightarrows M]$ is canonically
isomorphic to $\X$.
In summary, giving a groupoid space $X_1 \rightrightarrows X_0$
is equivalent to giving a Deligne-Mumford stack $\X$ and an atlas
$X_0 \to \X$.

Let $\xi:U \to \X$ be a morphism, which is considered as an object of $\X$.
If $\xi$ lifts to $\xi':U \to X_0$, then the automorphism group of $\xi$ is
identified with the set of morphisms $\eta: U \to X_1$ with 
$q_1 \circ \eta = q_2 \circ \eta = \xi'$.

\begin{defn}\label{def-morp-gropoidspace}
A \textit{morphism}
 $f:(Y_1 \rightrightarrows Y_0) \to (X_1 \rightrightarrows X_0)$
\textit{of groupoid spaces}
 is a pair of morphisms $f_i:Y_i \to X_i\; (i=0,1)$ which
respects the groupoid space structures.
\end{defn}

Given a morphism 
$f:(Y_1 \rightrightarrows Y_0) \to (X_1 \rightrightarrows X_0)$,
then we have a natural morphism of prestacks
$[f]':[Y_1 \rightrightarrows Y_0]' \to [X_1 \rightrightarrows X_0]' $ and
hence a natural morphism of stacks
$[f]:[Y_1 \rightrightarrows Y_0] \to [X_1 \rightrightarrows X_0] $ 
from \cite[Lem.\ 3.2]{LMB}.

Conversely, consider a commutative diagram 
\[\xymatrix{
    Y_0 \ar[d] \ar[r]^{f_0} & X_0 \ar[d]\\
    \Y   \ar[r]_g & \X
}\]
such that $\X$, $\Y$ are Deligne-Mumford stacks and
the vertical arrows are atlases.
If we define $Y_1 := Y_0 \times_\Y Y_0$ and $X_1 := X_0 \times_\X X_0$,
and if $f_1 : Y_1 \to X_1$ is the natural morphism,
then the pair of $ ( f_0 , f_1 ) $ determines 
a morphism $f:(Y_1 \rightrightarrows Y_0) \to (X_1 \rightrightarrows X_0)$
 of groupoid spaces. Evidently, $[f]=g$.

\begin{expl}
Let  $\X = [Z/G]$ be a quotient stack with $G$ finite.
There is a canonical atlas $Z \to \X$.
Then the groupoid space $Z \times _\X Z \rightrightarrows Z$ is isomorphic to 
the groupoid space 
\[\xymatrix{
Z \times G \ar@<-1ex>[rr]_{ \pr_1 } \ar@<1ex>^{ G-\text{action} }[rr]
& & Z
}\]
whose origin, inverse and composition are induced by the group structure of
$G$.
\end{expl}

Let $\X$ be a locally integral
Deligne-Mumford stack,  
associated to a groupoid space $X_1 \rightrightarrows X_0$. 
Let $X_i^\nor$ be the normalization of $X_i$ respectively.
Then the lifts of the structure morphisms of $X_1 \rightrightarrows X_0$
induce a groupoid space
 $ X_1^\nor \rightrightarrows X_0^\nor $.

\begin{defn}\cite[Def.\ 1.18]{Vistoli}\label{normalization}
We define the \textit{normalization} $\X^\nor$ of $\X$ to be 
the stack associated to $X_1^\nor \rightrightarrows X_0^\nor$.
\end{defn}

It is easy to show the uniqueness and the universality of the normalization.

\subsection{Quasi-coherent sheaves}

\begin{defn}\label{def-sheaf}
A \textit{quasi-coherent sheaf} $\cF$ on a Deligne-Mumford stack $\X$ 
consists of the following data:
\begin{enumerate}
    \item for each \'etale morphism $U \to \X$ with a scheme $U$, a 
 quasi-coherent sheaf $\cF_U$ on $U$,
    \item for each diagram of \'etale morphisms
  \[
\xymatrix{
V \ar[rr]^{\varphi} \ar[dr] && U  \ar[dl] \\
  & \X & 
}
  \]
  with $V$, $U$ schemes, an isomorphism 
$\varTheta _\varphi : \cF_V \to \varphi^* \cF_U$ which satisfies
 the cocycle condition.
\end{enumerate}
\end{defn}

\begin{expl}
\begin{enumerate}
	\item The \textit{structure sheaf} $\cO_\X$ on $\X$ is defined by 
         $ (\cO_\X)_U := \cO_U$.
        \item The \textit{sheaf of differentials} $\Omega_{\X/S}$ is defined by
        $ (\Omega_{\X/S})_U := \Omega_{U/S} $ and by the canonical isomorphism.
        \item Let $ f: \Y \to \X$ be a morphism of Deligne-Mumford stacks.
         We define the \textit{sheaf of relative differentials} 
       $\Omega_{\Y/\X}$
      of $\Y$ over $\X$ to be the unique sheaf such that 
      for each commutative diagram
\[
\xymatrix{
 V \ar[r]^{\txt{\'etale}} \ar[d] & \Y \ar[d] \\
 U \ar[r]_{\txt{\'etale}} & \X ,
}
\] 
        $(\Omega_{\Y/\X})_V := \Omega_{V/U}$. 
  Then we have the following exact sequence
\[
    f^* \Omega_{\X/S} \to \Omega_{\Y/S} \to \Omega_{\Y/\X} \to 0 .
\]
\end{enumerate}
\end{expl}

\begin{defn}
In Definition \ref{def-sheaf}, if every $\cF_U$ is an $\cO_U$-algebra and
 every $\varTheta _\varphi $ is a homomorphism of $ \cO_V$-algebras, then
we say that $\cF$ is an $\cO_\X$-algebra.
\end{defn}

As the case of schemes, to an $\cO_\X$-algebra $\cF$,
 we can associate a representable and affine morphism 
$\cSpec \cF \to \X$. For details, see \cite[(14.2)]{LMB}.

\subsection{Inertia stacks}\label{subsec-inertia}

In this subsection, we study the inertia stack.
It is an algebro-geometric object corresponding to the \textit{twisted sector}
which was introduced by Kawasaki \cite{kawasaki1} and used by Chen and Ruan
to define the orbifold cohomology \cite{CR}.

\begin{defn}
For a Deligne-Mumford stack $\X$,
its \textit{inertia stack}, denoted $I(\X)$, is the stack defined as follows:
an object over $ U \in \schS $ 
is a pair $(\xi,\alpha)$ where $\xi\in \ob\X_U$ and $\alpha\in 
\Aut(\xi)$, and a morphism $(\xi,\alpha)\to(\eta,\beta)$ is a morphism
$\gamma:\xi\to\eta$ in $\X$ such that $\gamma\circ\alpha=\beta\circ\gamma$.
\end{defn}

There is a natural forgetting morphism $I(\X)\to \X$.
The forgetting morphism $I(\X)\to \X$ is isomorphic
to 
\begin{equation*}
\pr_1:\X\times_{\Delta,\X\times\X,\Delta}\X\to\X.
\end{equation*}
Hence if $\X$ is complete, then so is $I(\X)$. 

The following lemma may be well-known:
\begin{lem}\label{local-inertia}
Let $Z$ be a scheme with an action of a finite group $G$.
Then we have an isomorphism:
\begin{equation*}
I([Z/G])\cong \coprod_{g\in \Conj(G)}[Z^g/C(g)],
\end{equation*}
where $\Conj(G)$ is a set of representatives of the conjugacy classes,
$Z^g$ is the locus of fixed points under the $g$-action and
$C(g)$ is the centralizer of $g$.
\end{lem}

\begin{proof}
Let $U\in\schS$ be a connected scheme.
An object of $[Z/G]$ over $U$ is a $G$-torsor $P\to U$ with a
$G$-equivariant morphism $P\to Z$. Its automorphism $\alpha$ is an
automorphism of a $G$-torsor $P\to U$ compatible with $P\to Z$.
For some \'{e}tale surjective $V\to U$, the fiber product 
$P_V:=P\times_U V$ is isomorphic to $G\times V$ as $G$-torsors over $V$.
Here $G\times V$ is a $G$-torsor for the \textit{right action} of $G$.
The pull back $\alpha_V$ is represented by the \textit{left action} of some
$g\in G$. $g$ is determined  up to conjugacy and we can assume 
$g\in\Conj(G)$.

Let $\psi:G\times V \cong P_V \to P$ be the natural morphism. Now let us show
that if $a\in C(g)$ and $b\notin C(g)$,
then we have $\psi(\{a\}\times V)\cap \psi(\{b\}\times V) =\emptyset $.
Let $x$ (resp.\ $y$) be  a geometric point of $\{a\}\times V$ 
(resp.\ $\{b\}\times V$) and assume $\psi(x)=\psi(y)$. Then we have
\begin{equation*}
\psi(x)=\psi(gxg^{-1})=\alpha \psi(x) g^{-1}=\alpha \psi(y) g^{-1}
=\psi(gyg^{-1})\ne\psi(y) .
\end{equation*}
This is a contradiction. So $P$ decomposes into $C(g)$-torsors
 as $P\cong P'\times J$ where 
$P':=\psi(C(g)\times V)$ and $J$ is a finite set.
Let $f_V$ denote the composition $f \circ \psi$.
In the following diagram,
\[ \xymatrix{
C(g)\times V \ar[r]^>>>>>{f_V} \ar[d]_{\alpha_V} & Z \ar[d]^g \\
C(g)\times V \ar[ur]_{f_V} \ar[r]_>>>>>{f_V} & Z
}\]
we have $g\circ f_V =f_V \circ \alpha_V$, since $\alpha_V$ equals the right
action of $g$ on $C(g)\times V$ and $f_V$ is $G$-equivariant. We also have
$f_V=f_V\circ\alpha_V$ and hence $g\circ f_V=f_V$. It implies that 
$f_V(C(g)\times V)$
is in $Z^g$ and hence so is $f(P')$. 
Thus the $C(g)$-torsor $P'\to V$ with $f:P'\to Z^g$ is an object
of $[Z^g/C(g)]$. For a non-connected $U\in\schS$ and an object of $I([Z/G])$
over $U$, we can assign it an object of $\coprod_{g\in\Conj(G)}[Z^g/C(g)]$,
in the obvious way. We leave the rest for the reader.
\end{proof}

\begin{defn}
A morphism $f:\Y \to \X$ of
stacks is \textit{barely faithful} if 
for every object $\xi$ of $\Y$, 
 the map $\Aut(\xi)\to\Aut(f(\xi))$ is bijective.
\end{defn}

Clearly, all barely faithful morphisms 
are  faithful functors.
From \cite[Prop.\ 2.3 and Cor.\ 8.1.2]{LMB}, all barely faithful morphisms
of Deligne-Mumford stacks
are \textit{representable} in the sense of \cite[D\'ef.\ 3.9]{LMB}.
Because all separated and quasi-finite morphisms of algebraic spaces are
sch\'ematique \cite[Th.\ A.2]{LMB}, 
all barely faithful and quasi-finite morphisms of 
Deligne-Mumford stacks are \textit{representable} for our definition.

\begin{expl}
All immersions are barely faithful. 
All morphisms of schemes are barely faithful.
\end{expl}

\begin{lem}\label{barely-base-change}
Barely faithful morphisms are stable under base change.
\end{lem}

\begin{proof}
Let $f:\Y \to \X$ be a barely representable morphism of stacks and 
$a:\X' \to \X$ any morphism of stacks.
An object of the fiber product $\Y \times_\X \X'$ is a triple $(\xi,\eta,
\alpha)$, where $\xi$ is an object of $\Y$, $\eta$ is an object of $\X'$ and
$\alpha : f(\xi) \to a(\eta)$ is a morphism in $\X_U$ for some $U\in\schS$.
Its automorphism is a pair of automorphisms $\varphi\in\Aut(\xi)$ and 
$\psi\in\Aut(\eta)$ with $\alpha\circ f(\varphi) = a(\psi)\circ\alpha$.
Since the map $f:\Aut(\xi) \to \Aut(f(\xi))$ is bijective, for each $\psi$,
there is one and only one
 $\varphi$ with  $\alpha\circ f(\varphi) = a(\psi)\circ\alpha$.
We have thus proved the lemma.
\end{proof}

\begin{prop}\label{barely-inertia}
Let $f:\Y \to \X$ be a barely faithful morphism of Deligne-Mumford stacks. 
Then the inertia stack $I(\Y)$ is naturally isomorphic to 
the fiber product $\Y \times_\X I(\X) $.
\end{prop}

\begin{proof}
The natural morphism $\Psi:I(\Y) \to \Y \times_\X I(\X) $ is 
defined as follows: for an object $\xi$ of $\Y$ and its automorphism $\alpha$,
the pair $(\xi , \alpha)$, which is an object of $I(\Y)$, is mapped to
the triple $(\xi,(f(\xi), f(\alpha)), \id_{f(\xi)} )$.

We will first show $\Psi$ is a fully faithful functor.
Let $\xi$ be an object of $\Y$.
The automorphism group of $(\xi,(f(\xi), f(\alpha)), \id_{f(\xi)} )$ is 
a pair of automorphisms 
$\beta \in \Aut(\xi)$ and 
$\nu \in \Aut((f(\xi), f(\alpha))) = C(f(\alpha))$ such that
$f(\beta) = \nu $. Hence $\Psi$ is barely faithful.
Let $\eta$ be another object of $\Y$ and $\tau$ an automorphism of $\eta$.
It suffices to show that  if 
$\Hom_{I(\Y)}((\xi,\alpha),(\eta,\tau))=\emptyset$,
then $\Hom_{\Y \times_\X I(\X)}(\Psi(\xi,\alpha),\Psi(\eta,\tau))=\emptyset$.
Suppose that  there is an element $(\sigma, \varepsilon)$ of 
$\Hom_{\Y \times_\X I(\X)}(\Psi(\xi,\alpha),\Psi(\eta,\tau))$,
where $\sigma$ is a morphism $\xi \to \eta$ and $\varepsilon$ is a morphism
$f(\xi)\to f(\eta)$ such that the diagram
\[\xymatrix{
 f(\xi)\ar[d]_{f(\alpha)} \ar[r]_{\varepsilon} & f(\eta) \ar[d]^{f(\tau)} \\
 f(\xi)\ar[r]_{\varepsilon} & f(\eta)
}\]
is commutative and $f(\sigma) = \varepsilon$.
Since $f$ is barely faithful, the diagram
\[\xymatrix{
 {\xi}\ar[d]_\alpha \ar[r]_{\sigma} & {\eta} \ar[d]^\tau \\
 {\xi}\ar[r]_{\sigma} & {\eta}
}\]
is commutative, that is, $\Hom_{I(\Y)}((\xi,\alpha),(\eta,\tau))\ne\emptyset$.

Now, let us show $\Psi$ is an isomorphism, that is, an equivalence of 
categories. Let $(\xi,(\theta, \beta), \upsilon )$ be an object of 
$\Y \times_\X I(\X)$ where $\upsilon$ is an isomorphism $f(\xi) \to \theta$.
Then there is a natural bijection $\Aut(\xi)\to \Aut(f(\xi))\to \Aut(\theta)$.
Let $\alpha\in\Aut(\xi)$ be an automorphism corresponding to 
$\beta\in\Aut(\theta)$.
Then we can see that $\Phi(\xi,\alpha)$ is isomorphic to 
$(\xi,(\theta, \beta), \upsilon )$. We have thus completed the proof.
\end{proof}

\begin{cor}\label{cor-inertia-smooth}
If $S= \Spec \CC$ and 
 $\X$ is a smooth Deligne-Mumford stack,
then  $I(\X)$ is also smooth.
\end{cor}

\begin{proof}
From Lemma \ref{locallyquotient} and Lemma \ref{barely-base-change},
 there is an \'etale, surjective
and barely faithful morphism
$\coprod_i [M_i/G_i] \to \X$
 such that each $M_i$ is smooth and each $G_i$ is a finite group. 
Then the assertion follows from Lemma \ref{local-inertia} and
Proposition \ref{barely-inertia}.
\end{proof}


\subsection{Coarse moduli space}\label{Coarse moduli space}

\begin{defn}
Let $\X$ be a Deligne-Mumford stack. The \textit{coarse moduli space}
of $\X$ is an algebraic space $X$ with a morphism $\X\to X$ such that:
\begin{enumerate}
	   \item for any algebraically closed field $\Omega$ with a morphism
$\Spec \Omega \to S$,  $\X(\Omega)\to X(\Omega)$ is a bijection,
        \item for any algebraic space
         $Y$, any morphism $\X\to Y$ uniquely
           factors as $\X \to X \to Y$.
\end{enumerate}
\end{defn}

Keel and Mori  proved that
the coarse moduli space always exists \cite[Cor.\ 1.3]{Keel-Mori}.

\begin{expl}
Let $Z$ be an algebraic space and $G$ a finite group acting on $Z$.
Then the coarse moduli space of the quotient stack $[Z/G]$ is the quotient
algebraic space $Z/G$. 
\end{expl}

The following lemma is well-known.

\begin{lem}\textit{e.g.} \cite[Lem.\ 2.2.3]{AV}\label{locallyquotient}
Let $\X$ be a Deligne-Mumford stack and $X$ its coarse moduli space.
Then there is an \'etale covering $(X_i \to X)_i$ such that 
$\X \times_X X_i$ is isomorphic to a quotient stack $[Z_i/G_i]$ with 
a scheme $Z_i$ and a finite group $G_i$. 
Hence the canonical morphism $\X \to X$ is proper.
\end{lem}

Now we assume $S = \Spec \CC$.

\begin{defn}\label{quot-sing}
Let $X$ be a variety.
We say that $X$ has \textit{quotient singularities} if there is an \'{e}tale
covering $(U_i/G_i \to X)_i$ with a smooth variety $U_i$ and a finite group
$G_i$.
\end{defn}

Lemma \ref{locallyquotient} shows that for 
a variety $X$,
$X$ has quotient singularities if it is the coarse moduli space of
some smooth Deligne-Mumford stack.
In fact, 
`only if' also holds (Lemma \ref{smoothing}).

Let $\X$ be a smooth Deligne-Mumford stack and
$x : \Spec \CC \to \X$  a closed point.
Then $\Aut (x)$ acts on the tangent space $T_x \X$.

\begin{defn}
We say that $\alpha \in \Aut (x)$ is a \textit{reflection} if 
the subspace of the $\alpha$-fixed points $(T_{x} \X)^{\alpha}$ is 
of codimension 1.
\end{defn}

\begin{lem}\label{smoothing}
Let $X$ be a $k$-variety with quotient singularities.
Then there is a smooth Deligne-Mumford stack $\X$ 
without reflections such that the automorphism
group of general geometric points is trivial and $X$ is the coarse moduli 
space of $\X$. 
\end{lem}

\begin{proof}
We will give only a sketch. There is a finite set of pairs $(V_i\to X, G_i)_i$
such that:
\begin{enumerate}
	\item $V_i$ is a smooth variety,
        \item $G_i$ is a finite group acting effectively on $V_i$ without
                   reflections,
        \item $V_i \to X$
         is a morphism \'{e}tale in codimension 1 which factors
                as $V_i\to V_i/G_i \to X$ with $V_i/G_i \to X$ \'{e}tale.
\end{enumerate}
Let $V_{ij}$ be the normalization of $V_i\times_X V_j$.
Then the natural morphisms $V_{ij} \to V_i$ and $V_{ij} \to V_j$ are
 \'{e}tale in codimension 1. 
From the purity of branch locus, they are actually \'{e}tale,
and hence $V_{ij}$ is smooth. 
The diagonal $\Delta: V_i\to V_i
\times_X V_i$ factors through $\Delta':V_i \to V_{ij}$. Then, with the
suitable multiplication morphism, the diagram
\[ \xymatrix{
{\coprod} V_{i,j} \ar@<-1ex>[r] \ar[r] & {\coprod} V_i \ar@<-1ex>[l]_{\Delta'} 
} \]
has the structure of groupoid space. We set $\X$ to be the associated stack.
Clearly $\X$ has no reflections.
The canonical morphism $\X\to X$ makes $X$ the coarse moduli space of $\X$
(see \cite[Prop.\ 9.2]{Gillet}). Any geometric point $x$ of $\X$ has a lift 
$\tilde x:\Spec \Omega \to V_i$ with $\Omega$ algebraically closed field. 
The automorphism group of $x$ is
identified with $\{y:\Spec\Omega \to V_{ii}\, 
              |\ p_1\circ y=p_2 \circ y=\tilde X
  \}$. If $y$ is over the smooth locus of $X$, then this group is trivial.
\end{proof}


\def\cprime{$'$}

\end{document}